\documentclass{article}

\usepackage[utf8]{inputenc}
\usepackage[letterpaper,top=2cm,bottom=2cm,left=3cm,right=3cm,marginparwidth=1.75cm]{geometry}

\usepackage{amsmath,amsfonts,amsthm}
\usepackage{cleveref}
\usepackage{mathtools}
\usepackage{tikz}
\usepackage{makecell}

\newcommand\R{\mathbb{R}}
\newcommand\N{\mathbb{N}}

\newcommand{\rank}{\mathrm{rank}}
\renewcommand{\P}{\mathbb{P}}
\newcommand{\M}{\mathcal{M}}

\newtheorem{theorem}{Theorem}
\newtheorem{lemma}[theorem]{Lemma}

\newtheorem{proposition}[theorem]{Proposition}
\newtheorem{corollary}[theorem]{Corollary}
\theoremstyle{definition}
\newtheorem{definition}[theorem]{Definition}
\newtheorem{remark}[theorem]{Remark}

\newtheorem{example}[theorem]{Example}

\def\dist{\operatorname{dist}}

\title{A semialgebraic description of matrices of positive semidefinite rank at most $2$}
\author{Kaie Kubjas, Lilja Metsälampi, and Milo Orlich}

\begin{document}
\maketitle

\begin{abstract}
    The main result of this paper is a semialgebraic description of the set of matrices of rank at most $3$ and psd rank at most $2$. Such a semialgebraic description allows one to check whether a matrix has rank at most $3$ and psd rank at most $2$. This description builds on the geometric characterization of matrices of rank $3$ and psd rank $2$ using the existence of an ellipse nested between two polygons. The main step towards deriving the semialgebraic description is the observation that we need to consider only finitely many ellipses that go through some of the vertices of the inner polygon and are tangent to some of the edges of the outer polygon. Although the description is not completely quantifier-free, for a general matrix, the quantified variables satisfy a system of polynomial equations that has at most four solutions corresponding to reduced conics. 
\end{abstract}

\section{Introduction}

A size-$k$ positive semidefinite factorization of a nonnegative matrix $M=(m_{ij}) \in \R_{\geq 0}^{p \times q}$ is a collection $(A_1,\dots,A_p,\allowbreak B_1 \dots,B_q)$ of $k \times k$ symmetric positive semidefinite matrices such that $m_{ij} = \langle A_i, B_j \rangle$ where $\langle \cdot, \cdot \rangle$ denotes the trace inner product.
The \textit{positive semidefinite (psd) rank} of $M$ is the smallest $k \in \N$ such that $M$ admits a size-$k$ psd factorization. We denote it by $\operatorname{rank}_{psd}(M)$.

Psd rank was first introduced to study the semidefinite extension complexity of polytopes, a notion measuring the minimal complexity of solving a semidefinite program over a polytope~\cite{fiorini2012linear}. An information theoretic application of psd rank is given in~\cite{jain2013efficient}.
A comprehensive review of psd factorizations and psd rank can be found in \cite{fawzi2015positive}. Psd factorizations are closely related to nonnegative factorizations and they are both examples of a more general class of cone factorizations introduced in \cite{gouveia2013lifts}.
Work on psd factorizations spans several directions: computational complexity~\cite{shitov2017complexity,shitov2018matrices,shitov2024further}, algorithms for computing psd factorizations~\cite{lahat2021positive,soh2021non,vandaele2018algorithms}, uniqueness of psd factorizations~\cite{dawson2024uniquenesssize2positivesemidefinite}, and the psd rank of polytopes~\cite{gouveia2013polytopes,gouveia2015worst,gouveia2017four}. 

Let $\M_{r,k}^{p \times q}$ denote the set of matrices of rank at most $r$ and psd rank at most $k$. This set is semialgebraic by the Tarski--Seidenberg theorem. In this paper, we aim to characterize $\M_{r,k}^{p \times q}$ for $r=3$ and $k=2$. Such a description can be used to determine whether a nonnegative matrix of rank at most~$3$ has psd rank at most~$2$ without computing a size-$2$ psd factorization.
We start by justifying why we consider the setting $r=3$ and $k=2$.

\begin{proposition}[{{\cite[Proposition 2]{fawzi2015positive}}}]\label{lem:relations-between-ranks}
    If\/ $M \in \R_{\geq 0}^{p \times q}$ is a nonnegative matrix, then 
    \begin{equation*}\label{eq:relations-between-ranks}
        \frac{1}{2}\sqrt{1+8\operatorname{rank}(M)}-\frac{1}{2} \leq \operatorname{rank}_{psd}(M) \leq \operatorname{rank}_+(M) \leq \min(p,q),
    \end{equation*}
    where $\operatorname{rank}_+(M)$ denotes the nonnegative rank of $M$.
\end{proposition}

\begin{proposition}[{{\cite[Proposition 3]{fawzi2015positive}}}]\label{prop:psd_rank_for_matrices_of_rank_one_and_two}
If\/ $M \in \R_{\geq 0}^{p \times q}$ is a nonnegative matrix, then 
\[
\operatorname{rank}(M)=1 \Leftrightarrow \operatorname{rank}_{psd}(M)=1 \quad \quad \text{and} \quad \quad \operatorname{rank}(M)=2 \Rightarrow \operatorname{rank}_{psd}(M)=2.
\]
\end{proposition}

By~\Cref{prop:psd_rank_for_matrices_of_rank_one_and_two}, the set of matrices of rank at most $1$ (resp.~$2$) and psd rank at most $1$ (resp.~$2$) is characterized by the vanishing of $2 \times 2$ (resp.~$3 \times 3$) minors of the matrix and the entries of the matrix being nonnegative. 
\Cref{lem:relations-between-ranks} gives that for matrices of rank $3$ the psd rank is at least $2$. Therefore, the first nontrivial case is the set of matrices of rank at most $3$ and psd rank at most $2$. \Cref{lem:relations-between-ranks,prop:psd_rank_for_matrices_of_rank_one_and_two} also give that $\operatorname{rank}_{psd}(M)=2 \Rightarrow \operatorname{rank}(M) \in \{2,3\}$.  Therefore, a semialgebraic description for the set of matrices of rank at most $3$ and psd rank at most $2$ also gives a semialgebraic description for matrices of psd rank at most $2$ without a rank constraint.

Our main tool for characterizing $\M_{3,2}^{p \times q}$ is a geometric description of matrices of rank $3$ and psd rank $2$. This description associates to $M \in \R_{\geq 0}^{p \times q}$ of rank $3$ two nested polygons $P,Q$ and states that its psd rank is $2$ if and only if there exists an ellipse $E$ nested between the two polygons~\cite[Proposition 9]{gouveia2015worst}, see~\Cref{subsec:geometry-of-psd} for details. Our first main result that allows to derive a semialgebraic description is that instead of considering all possible ellipses in $\mathbb{R}^2$, we have to consider only finitely many ellipse candidates:

\begin{theorem}\label{thm:main}
Let $P, Q \subset \mathbb{R}^2$ be polygons and $E \subset \mathbb{R}^2$ an ellipse such that $P \subset \operatorname{conv}(E) \subset Q$. Then there exists an ellipse $E'$ such that $P \subset \operatorname{conv}(E') \subset Q$ and $E'$ has at least five incidence and tangency relations with the polygons~$P$ and~$Q$, respectively.
\end{theorem}

The inspiration for this theorem comes from~\cite{kubjas2018positive}, where the boundaries of the set $\M_{3,2}^{p \times q}$ are studied. In particular, it is shown in~\cite[Theorem 3.6]{kubjas2018positive} that a matrix $M$ is on the boundary of $\M_{3,2}^{p \times q}$ if and only if every ellipse that is nested between the associated polygons $P$ and $Q$ contains at least three vertices of the inner polygon $P$ and is tangent to at least three edges of the outer polygon~$Q$. The idea behind our theorem is that if $M$ is a matrix in $\M_{3,2}^{p \times q}$ but not necessarily on the boundary, then one less total incidence and tangency relation suffices. \Cref{thm:main} is used to derive our second main result, the semialgebraic description in~\Cref{thm:semialgebraic-description} (\Cref{sec:semialgebraic_description}).

The proof of~\Cref{thm:main} has two main steps. In~\Cref{prop:increasing_incidences} (\Cref{sec:families-of-ellipses}), we show that if a path-connected family of conics with fixed minimal incidence relations with $P$ and tangency relations with $Q$ contains an ellipse $E$ such that $P \subset \operatorname{conv}(E) \subset Q$ as well as a conic that is not an ellipse, then it also contains an ellipse $E'$ such that $P \subset \operatorname{conv}(E') \subset Q$ and $E'$ has an additional incidence with $P$ or tangency with $Q$. Then in the proof of~\Cref{thm:main} (\Cref{sec:infinite_to_finite}) we show that for any four or fewer incidence and tangency relations and every ellipse $E$ such that $P \subset \operatorname{conv}(E) \subset Q$ satisfying these relations, we can find a conic that is not an ellipse in the path-connected component of $E$ within the family of conics with at least these incidence and tangency relations. This amounts to studying the geometry of the family of conics with given incidence and tangency relations.

Related work on this subject includes \cite[Proposition 13]{gouveia2015worst} and \cite[Section 4]{fawzi2015positive}, where deciding whether a given matrix has psd rank $k$ and $2$, respectively, are given as semialgebraic set feasibility problems. In~\cite{shitov2018matrices}, the results of Renegar on the first-order
theory of the reals are used to give a first-order quantifier-free formula for the Nested Projection of Affine Slice of Cone problem and to prove that there exists a polynomial-time algorithm which decides whether a matrix has some fixed psd rank. Our semialgebraic description is of yet another flavor: it is not quantifier-free, but for each fixed matrix to find the quantified variables, it is enough to solve finitely many systems of polynomial equations and then check whether any of the solutions satisfy finitely many polynomial inequalities. These systems correspond to finding conics that go through $k$ vertices of $P$ and are tangent to $5-k$ edges of $Q$, and by~\cite{bashelor2008enumerative}, if the points and edges are in general position, each of these systems has at most four solutions corresponding to reduced conics. Our semialgebraic description has similarities with the semialgebraic description for matrices of nonnegative rank at most 3 in~\cite{kubjas2015fixed}.

The outline of the paper is the following. In~\Cref{sec:background}, we recall necessary background about psd rank and ellipses. In~\Cref{sec:families-of-ellipses}, we study properties of families of conics such that each conic is incident to given points and tangent to given lines. In~\Cref{sec:infinite_to_finite}, we show that to determine whether there exists an ellipse nested between two polygons, we can check this for finitely many ellipses. In~\Cref{sec:semialgebraic_description}, we use this to derive a semialgebraic description for matrices of rank at most~$3$ and psd rank at most~$2$.

\section{Background}\label{sec:background}

In this section we give the necessary background for deriving our results.
In \Cref{subsec:geometry-of-psd} we introduce a geometric interpretation of matrices of rank $3$ and psd rank $2$ using ellipses nested between polygons.
In \Cref{subsec:ellipses} we recall some basic results on conics and ellipses.

\subsection{Geometry of psd rank}\label{subsec:geometry-of-psd}

Let $M \in \R_{\geq 0}^{p \times q}$ be a nonnegative matrix of rank $r$.
In this subsection, we show how to associate to $M$ a pair $P \subseteq Q \subset \R^{r-1}$ where $P$ is a polytope and $Q$ is a polyhedron. Moreover, when $M$ has rank $3$, we recall a geometric characterization of $M$ being psd rank $2$.
For matrices $M$ satisfying $M\mathbf{1}=\mathbf{1}$, the polyhedron $Q$ can be chosen to be bounded as well by \cite[Theorem 1]{gillis2012geometric} and \cite[Lemma 4.1]{fawzi2015positive}. In fact, without loss of generality, we can consider matrices satisfying $M\mathbf{1}=\mathbf{1}$. Assume that $M$ has no zero rows.
Thus, the row sums of $M$ are positive, and, since scaling the rows of the matrix does not change the psd rank of the matrix by \cite[Theorem 1]{fawzi2015positive}, we can, without loss of generality, consider matrices $M$ satisfying $M\mathbf{1}= \mathbf{1}$.

\begin{definition}\label{def:nested-polytopes-of-M}
Let $M \in \R_{\geq 0}^{p \times q}$ be a nonnegative matrix of rank $r$ such that $M\mathbf{1}=\mathbf{1}$.
Consider a rank-size factorization $M = AB$ with $A$ having rows $A_i = (a_i,1)^T$. 
Let 
\begin{equation*}
    P = \operatorname{conv} (a_1, \dots , a_p)  \ \text{ and } \ Q = \{x \in \R^{r-1} \mid (x^T,1)B \geq 0 \}.
\end{equation*}
We call the polytopes $P$ and $Q$ the polytopes associated to the matrix $M$. 
\end{definition}

The psd rank of $M$ can be characterized in terms of $P$ and $Q$. The general setting is characterized in~\cite[Proposition 7]{gouveia2015worst}. We state the geometric characterization relevant for this paper, namely for checking whether a matrix of rank $3$ has psd rank $2$. In the rest of the paper, a polygon means a 2-dimensional polytope. In particular, a polygon in $\mathbb{R}^2$ has nonempty interior and it is bounded.

\begin{proposition} [{{\cite[Proposition 9]{gouveia2015worst}}}]\label{prop:psd-rank-2-iff-ellipse}
    Let\/ $M \in \R_{\geq 0}^{p \times q}$ be a nonnegative matrix of rank $3$ satisfying $M\mathbf{1}=\mathbf{1}$. Let $P \subseteq Q$ be the polygons associated to $M$ as in \Cref{def:nested-polytopes-of-M}.
    Then $M$ has psd rank $2$ if and only if there exists an ellipse $E \subset \R^2$ whose convex hull $\operatorname{conv}(E)$ satisfies $P \subset \operatorname{conv}(E) \subset Q$.
\end{proposition}

If there exists a zero entry in the matrix $M$, this means that a vertex of the polygon $P$ touches an edge of the polygon $Q$.
Hence, any ellipse nested between $P$ and $Q$ must pass through the vertex of $P$ and be tangent to the edge of $Q$.
If, however, all entries of $M$ are positive, then the boundaries of $P$ and $Q$ have no points in common.

\subsection{Conics}\label{subsec:ellipses}

\Cref{prop:psd-rank-2-iff-ellipse} will be the main tool in later sections for deriving our results. Motivated by this, we recall some results on conics and ellipses.

\begin{definition}\label{def:real-ellipse}
A \textit{conic} is the set of points in $\R^2$ that satisfy a quadratic equation in the variables $x,y$ of the form
$$u_0x^2+u_1xy+u_2y^2+u_3x+u_4y+u_5=0,$$
where not all of the coefficients $u_0,u_1,u_2$ are zero. We denote the polynomial on the left hand side by $U(x,y)$.
We define 
\[
S = \begin{pmatrix}
u_0 & u_1/2 & u_3/2\\
u_1/2 & u_2 & u_4/2\\
u_3/2 & u_4/2 & u_5
\end{pmatrix},
\quad \quad
S'=\begin{pmatrix}
         u_0 & u_1/2\\
         u_1/2 & u_2
\end{pmatrix}.
\]
The matrix $S$ is called \textit{matrix of the conic}.
The equation of a conic can be rewritten as 
\[
\begin{pmatrix}
x & y & 1
\end{pmatrix}
S
\begin{pmatrix}
x\\
y\\
1
\end{pmatrix}=0.
\]
\end{definition}

Conics are invariant under any nonzero scaling of $(u_0,\ldots,u_5) \in \mathbb{R}^6$, and hence can be considered as points $[u_0:\ldots:u_5]$ in the projective space $\mathbb{P}^5$, where at least one of $u_0,u_1,u_2$ is nonzero.

\begin{definition}
A \textit{real ellipse} (or simply \textit{ellipse} in the following, unless stated otherwise) is the locus of points in $\mathbb{R}^2$ such that the sum of distances from two fixed points in $\mathbb{R}^2$ is constant. The two fixed points are called the \textit{foci} or \textit{focal points} of the ellipse. Denoting the two focal points by $F_1$ and $F_2$, the line through $F_1$ and $F_2$ intersects the ellipse in two points, which are the extremes of a line segment that is called the \textit{major axis} of the ellipse. The line perpendicular to the major axis that goes through the midpoint between $F_1$ and $F_2$ also cuts the ellipse in two points, which are the extremes of a segment called the \textit{minor axis} of the ellipse.

In the special case that the ellipse is a circle, so that the focal points are the same point, the major and minor axes are not defined, but we shall nonetheless talk about their lengths, with a little abuse of notation, by meaning the diameter of the circle.
\end{definition}

\begin{lemma}\label{lem:coeffs-of-ellipse}
Let $F_1=(a,b)$ and $F_2=(c,d)$ be points, not necessarily distinct. 
Consider the ellipse that consists of all points such that the sum of distances from the focal points $F_1$ and $F_2$ equals $e\ge0$. 
If this ellipse is represented in the form
\begin{equation}\label{eq:ellipse-as-usual-conic}
u_0x^2+u_1xy+u_2y^2+u_3x+u_4y+u_5=0,
\end{equation}
then the coefficients $u_i$ are (up to a scalar multiple)
\begin{equation*}\label{eq:u_coefficients_in_terms_of_foci_and_distance}
    \begin{split}
        u_0&=e^2-(a-c)^2,\\
        u_1&=-2(a-c)(b-d),\\
        u_2&=e^2-(b-d)^2,\\
        u_3&=(a-c)(b^2-d^2)-(a+c)(e^2-(a-c)^2)=-\frac{1}{2}u_1 (b+d)-(a+c)u_0,\\
        u_4&=(b - d)(a^2 - c^2)-(b + d) (e^2-(b - d)^2)=-\frac{1}{2}u_1 (a+c)-(b + d) u_2,\\
        u_5 &= - 1/4 (e^2-a^2 - b^2 - c^2 - d^2)^2+(a^2+b^2)(c^2+d^2),\\
            &= - 1/4 (e^2-\|F_1\|^2-\|F_2\|^2)^2+\|F_1\|^2\|F_2\|^2,
    \end{split}
\end{equation*}
where in the last line we denote $\|F_1\|^2:=a^2+b^2$ and $\|F_2\|^2:=c^2+d^2$.
\end{lemma}

\begin{proof}
By definition, the ellipse in question is the locus of points $(p_1,p_2) \in \mathbb{R}^2$ that satisfy the equation
$$\sqrt{(p_1-a)^2+(p_2-b)^2}+\sqrt{(p_1-c)^2+(p_2-d)^2}=e.$$
By clearing out the square roots and rearranging the terms, a straightforward computation shows that  the ellipse  can be written in the form \eqref{eq:ellipse-as-usual-conic}, where the coefficients are as in the statement.
\end{proof}

\begin{remark}\label{rmk:dist-and-foci-assumption}
In the above lemma we implicitly assume that the number $e$ is strictly larger than the distance between $F_1$ and $F_2$, so that we have a real ellipse, i.e.,
\begin{equation}\label{eq:ebiggerthandist}
e>\dist(F_1,F_2):=\sqrt{(a-c)^2+(b-d)^2}.
\end{equation}
This is the only condition to get a real ellipse, as $a$, $b$, $c$ and $d$ can be any real numbers.
\end{remark}

\begin{corollary}\label{cor:point-inside-ellipse}
Consider an ellipse defined by\/ $U(x,y)$, where $u_0>0$.  The point $(p_1,p_2) \in \mathbb{R}^2$ is contained in the convex hull of the ellipse if and only if\/ $U(p_1,p_2) \leq 0$.
\end{corollary}

\begin{proof}
The quadratic part of the function $U(x,y)$ is strictly convex because  $u_0>0$ and $\det(S')>0$ imply that $u_2>0$, and hence its Hessian is positive definite. The function $U(x,y)$ is strictly convex because it is the sum of a strictly convex function and a linear function. Therefore, the convex hull of the ellipse defined by $U(x,y)$ is $\{(p_1,p_2) \in \mathbb{R}^2: U(p_1,p_2) \leq 0)\}$.
\end{proof}

\begin{remark}\label{rmk:real_ellipse_nonzero_coefficients}
    From~\eqref{eq:ebiggerthandist} and the expressions for $u_0$ and $u_2$ in~\Cref{lem:coeffs-of-ellipse} it follows that for a real ellipse $u_0 \neq 0$ and $u_2 \neq 0$. Moreover, from the characterization of ellipses in terms of $S$ and $S'$ below, the sign of $u_0$ has to be equal to the sign of $u_2$.
\end{remark}

The characterization of conics in terms of the associated matrices $S$ and $S'$ is given in \Cref{tab:classification_of_conics} and can be found in for example~\cite[p. 70]{spain2007analytical}. In particular, the conic with matrix $S$ is a real ellipse if and only if
    \begin{equation}\label{eqn:real-ellipse-conditions}
        \det(S') >0 \quad \text{ and } \quad (u_0+u_2)\det(S)<0.
    \end{equation}
A conic is called degenerate if and only if $\det(S)=0$.
\begin{table}[ht]
    \centering
    \begin{tabular}{|c|c|c|c|c|}
        \hline
         & $\det (S)$ & $\det (S')$ & $(u_0+u_2)\det S$ & $u_3^2+u_4^2-4(u_0+u_2)u_5$\\
        \hline
        \makecell{real ellipse} & $\neq 0$ & $>0$ & $<0$ & \\
        \hline
         \makecell{imaginary ellipse} & $\neq 0$ & $>0$ & $> 0$ & \\
        \hline
        \makecell{parabola} & $\neq 0$ & $=0$ &  &\\
        \hline
        \makecell{hyperbola} & $\neq 0$ & $<0$ &  &\\
        \hline
        \makecell{real intersecting lines} & $=0$ & $<0$  &  &\\
        \hline
        \makecell{conjugate complex intersecting lines} & $ = 0$ & $>0$ &  &\\
        \hline
        \makecell{real distinct parallel lines} & $= 0$ & $=0$ &  & $> 0$\\
        \hline
        \makecell{conjugate complex parallel lines} & $= 0$ & $=0$ &  & $< 0$\\
        \hline
        \makecell{coincident lines} & $= 0$ & $=0$ &  & $= 0$\\
        \hline
    \end{tabular}
    \caption{Classification of conics.}
    \label{tab:classification_of_conics}
\end{table}

We finish this section by characterizing when a conic is tangent to a line and when a conic intersects a line. We will consider tangency and intersections in the projective plane, to allow nice algebraic characterizations.  This does not change anything for real ellipses as even when considered in the projective plane, none of the intersection or tangency points are at infinity. We will consider tangency in the algebraic sense that substituting the equation of a line into the equation of the conic gives a quadratic equation with a double solution. For example, with this description two intersecting lines are tangent to any line that goes through the intersection point. 

\begin{definition}
Consider a conic and line defined respectively by the polynomials
\begin{equation*}
U(x,y)=u_0x^2+u_1xy+u_2y^2+u_3x+u_4y+u_5\qquad\text{and}\qquad
\ell(x,y)=rx+sy+t.
\end{equation*}
The \textit{tact invariant of\/ $U$ and\/ $\ell$} is defined as
\begin{equation*}
    \begin{split}
    T(U,\ell):=\ &
    - 4u_0u_2t^2
    + 4u_0u_4st 
    - 4u_0u_5s^2 
    - 2u_1u_3st 
    - 2u_1u_4rt
    + 4u_1u_5rs\\
    & + u_1^2t^2 
    + 4u_2u_3rt 
    - 4u_2u_5r^2 
    - 2u_3u_4rs 
    + u_3^2s^2 
    + u_4^2r^2.
    \end{split}
\end{equation*}
\end{definition}

\begin{lemma}\label{lem:inside-intersecting-ellipse}
The projectivizations of the conic and the line defined by\/ $U(x,y)$ and\/ $\ell(x,y)$ respectively intersect if and only if\/ $T(U,\ell)\ge0$. They are tangent if and only if\/ $T(U,\ell) = 0$.
\end{lemma}

\begin{proof}
The projectivizations of the conic and line are defined by 
\begin{equation*}
U(x,y,z)=u_0x^2+u_1xy+u_2y^2+u_3xz+u_4yz+u_5z^2\qquad\text{and}\qquad
\ell(x,y,z)=rx+sy+tz.
\end{equation*}
Since we started with a line in $\mathbb{R}^2$, at least one of $r$ and $s$ has to be non-zero, say $r$. If we substitute $x$ from the equation of the line to the equation of the conic, we get the single second-degree equation 
         $$(u_0s^2-u_1rs+u_2r^2)y^2+(2u_0st-u_1rt-u_3rs+u_4r^2)yz+(u_0t^2-u_3rt+u_5r^2)z^2=0.$$ 
         In order for this equation to have at least one real solution in $\mathbb{P}^1$, its discriminant has to be $\ge0$, and there is a double solution when the discriminant is exactly zero. A straightforward computation shows that the discriminant is equal to the tact invariant up to the positive factor $r^2$. This is done assuming that $r\ne0$, and analogously we could eliminate $y$ by assuming that $s\ne0$.
\end{proof}

An incidence relation between a conic and the polygon $P$ means that the ellipse passes through a vertex of $P$.
A tangency relation between a conic and $Q$ always means that the conic is tangent to an edge of $Q$.
When we say that a conic has five incidence and tangency relations with $P$ and $Q$ we mean that the ellipse has $k$ incidences with $P$ and $5-k$ tangencies with $Q$ for $0 \leq k \leq 5$.

\section{Families of conics with given incidence and tangency conditions}\label{sec:families-of-ellipses}

The main goal of this section is to prove~\Cref{prop:increasing_incidences}, which will be used in~\Cref{sec:infinite_to_finite} to show that if there is an ellipse nested between polygons $P$ and $Q$, then there is an ellipse with at least five incidence and tangency relations with $P$ and $Q$ nested between these polygons. The proof of~\Cref{prop:increasing_incidences} builds on~\Cref{cor:path_to_different_conic_implies_non_nestedness} and~\Cref{proposition:one_more_incidence}.

\begin{proposition}\label{prop:increasing_incidences}
Let $P \subset Q$ be polygons. Let\/ $\mathcal{C}$ be a path-connected family of conics such that each conic in the family contains  $p$ given vertices of\/ $P$ and is tangent to  $q$ given edges of\/ $Q$. If\/ $\mathcal{C}$ contains an ellipse $E$ such that $P \subset \operatorname{conv}(E) \subset Q$ and $E$ has no further incidence and tangency relations with $P$ and $Q$, and a conic that is not an ellipse, then it also contains an ellipse $E'$ such that $P \subset \operatorname{conv}(E') \subset Q$ and $E'$ has an additional incidence or tangency relation with $P$ or $Q$.
\end{proposition}

\begin{proposition}\label{prop:limits-of-ellipses}
    Let $C_1 \subset C_2$ be two fixed nested circles. Let $\mathcal C \subset \mathbb{P}^5$ be the space of conics  and let $\pi:[0,1]\to \mathcal C$ be a continuous path. 
    If $\pi$ is such that 
    \begin{itemize}
        \item $\pi(t)$ is an ellipse for all\/ $t\in[0,1)$, and
        \item $\pi(1)$ is a degenerate conic or a parabola,
    \end{itemize}
then there exists $t \in [0,1)$ such that the ellipse $\pi(t)$ is not nested between $C_1$ and $C_2$.
\end{proposition}

\begin{proof}
     We show the contrapositive of the statement, that is, if for every $t \in [0,1)$ the ellipse $\pi(t)$  is nested between $C_1$ and $C_2$, then $\pi(1)$ cannot be a parabola or a degenerate conic.
    Since $\pi(t)$ is an ellipse for $t \in [0,1)$, it can be lifted to a point $u(t)= (u_0(t), \dots, u_5(t)) \in \R^6\setminus\{0\}$ where $u_0, \dots, u_5$ are as in~\Cref{lem:coeffs-of-ellipse}.
    Let $s_1(t) \geq s_2(t)$ be the semi-axes of the ellipse $\pi(t)$.
    Then, the following identities hold:
    \begin{equation*}
        s_1(t) = \frac{1}{2}e(t) \text{ and } s_2(t) = \sqrt{\frac{1}{4}(e(t)^2-\operatorname{dist}(F_1(t),F_2(t))^2)}.
    \end{equation*}

    Assume that for every $t \in [0,1)$ the ellipse $\pi(t)$ is nested between the circles $C_1$ and $C_2$.
    Then there exists some $\alpha, \beta \in (0,\infty)$ such that $\beta \leq  s_2(t) \leq s_1(t) \leq \alpha$.
    Using the coefficients from~\Cref{lem:coeffs-of-ellipse}, for every $t \in[0,1)$
    \begin{equation*}
            \det(S_t') = e(t)^2(e(t)^2-\operatorname{dist}(F_1(t),F_2(t))^2), \text{ and } \det(S_t) = -\frac{1}{4}\det(S_t')^2.
    \end{equation*}
    Since $e(t) = 2 s_1(t)$ and $2s_2(t) = \sqrt{e(t)^2-\operatorname{dist}(F_1,F_2)^2}$, it follows that $\det(S_t') = 16s_1(t)^2s_2(t)^2$ and thus $\det(S_t) = -64 s_1(t)^4 s_2(t)^4$.
    Therefore, $\det(S_t') \in [16 \beta^4 , 16 \alpha^4]$ and $\det(S_t) \in [-64 \alpha^8, -64 \beta^8]$ for all $t \in[0,1)$.
    Moreover, $0$ is not contained in either of these intervals.

    Next, we show that the squared Euclidean norms $\|u(t)\|_2^2 = u_0(t)^2 + \cdots + u_5(t)^2$ have an upper and lower bound for each $t \in [0,1)$.
    Since $a(t), b(t),c(t),d(t)$ are the coordinates of the two focal points $F_1(t)$ and $F_2(t)$ which live inside a circle, it follows that for all $t \in [0,1)$ they are bounded.
    Similarly, since $e(t) = 2s_1(t)$ it follows that $2\beta \leq e(t) \leq 2 \alpha$.
    Therefore, each of the $|u_0(t)|, \dots, |u_5(t)|$ has an upper bound, so the squared norm $\|u(t)\|_2^2 \leq L^2$ for some $L > 0$.
    To derive the lower bound, we consider the sum $u_0(t)^2+ u_1(t)^2+u_2(t)^2$.
    A direct computation shows that
    \begin{equation*}
       \begin{split}
           u_0(t)^2+ u_1(t)^2+u_2(t)^2 
           = &(\operatorname{dist}(F_1(t),F_2(t))^2+ 4s_2^2)^2 + 2 (a(t)-c(t))^2(b(t)-d(t))^2 + 16 s_2(t)^4 \\
           \geq& 16 s_2(t)^4.
       \end{split}
    \end{equation*}
    Therefore, $\|u(t)\|_2^2 \geq u_0(t)^2 + u_1(t)^2 +u_2(t)^2 \geq 16s_2^4 \geq 16 \beta^4$.

    We now return to consider the original path $\pi$.
    Since the unit sphere $S^5 \subset \mathbb{R}^6$ is a covering space of $\P^5$, the path $\pi$ lifts uniquely to a path $\tilde \pi \colon [0,1] \to S^5$ such that $\tilde \pi(0)$ is a lift of $\pi(0)$ to the sphere with $u_0 > 0$. 
    Then $\tilde \pi(0)= \frac{u(0)}{\|u(0)\|_2}$, and indeed, for all $t \in [0,1)$ we have $\tilde \pi(t) = \frac{u(t)}{\|u(t)\|_2}$. Let $\det(\tilde S_t)$ and $\det(\tilde S_t')$ be the determinants associated to $\tilde \pi(t)$.
    For $t \in [0,1)$ we have $\det(\tilde S_t) = \frac{\det (S_t)}{\|u(t)\|_2^3}$ and $\det(\tilde S_t') = \frac{\det (S_t')}{\|u(t)\|_2^2}$.
    We have derived the following bounds:
    \begin{equation*}
        \begin{split}
            &\frac{-\alpha^8}{\beta^6} \leq \det(\tilde S_t) \leq -\frac{64 \beta^8}{L^3} < 0\\
            &0 < \frac{16 \beta^4}{L^2} \leq \det(\tilde S_t') \leq \frac{\alpha^4}{ \beta^4}.
        \end{split}
    \end{equation*}

    With the coordinates normalized to the sphere $S^5$, the determinants $\det(\tilde S_t)$ and $\det(\tilde S_t')$ are well-defined and continuous functions.
    In particular, $\lim_{t \to 1} \det(\tilde S_t) = \det(\tilde S_1)$ and $\lim_{t \to 1} \det(\tilde S_t') = \det(\tilde S_1')$.
    Since for $t \in [0,1)$ the values of $\det(\tilde S_t)$ and $\det(\tilde S_t')$ are contained in closed intervals not containing zero, it must be that $\lim_{t \to 1} \det(\tilde S_t)$ and $\lim_{t \to 1} \det(\tilde S_t')$ are contained in these intervals as well.
    In particular, $\det(\tilde S_1) \neq 0$ and $\det(\tilde S_1') \neq 0$, that is, $\pi(1)$ cannot be a parabola or a degenerate conic.
\end{proof}

\begin{corollary} \label{cor:path_to_different_conic_implies_non_nestedness}
    Let $\mathcal C$ be the space of conics and $P\subset Q$ be nested polygons in $\R^2$. 
    Let $\pi:[0,1]\to \mathcal C$ be a path in the space of conics, and assume that 
    \begin{itemize}
        \item $\pi(0)$ is an ellipse nested between $P$ and\/ $Q$, and 
        \item $\pi(1)$ is not an ellipse.
    \end{itemize}
    Then there is some $t_0\in(0,1)$ such that $\pi(t)$ is an ellipse for all\/ $t \in [0,t_0]$ and $\pi(t_0)$ is not nested between $P$ and\/ $Q$. 
\end{corollary}

\begin{proof}
    
Let
\[
T = \left\{ t \in [0,1]: \pi(t) \text{ is not an ellipse}\right\}.
\]
This set is closed in $\mathbb{R}$ with respect to the Euclidean topology, since not being an ellipse can be characterized using inequalities $\det(S') \leq 0$ or $(u_0+u_2) \det S\geq 0$. This set is nonempty as $1 \in T$. Every nonempty, closed and bounded set has a minimum element. Let $t^*$ be this minimum. We claim that $\pi(t^*)$ is either a degenerate conic or a parabola. Indeed, for each $t$  let us fix a representative of $\pi(t)$ of norm~$1$, call it $\tilde\pi(t)$, in such a way that $\tilde\pi:[0,1]\to\R^6$ is continuous.

    The path $\pi$ cannot go from a real ellipse to a hyperbola or an imaginary ellipse without going through a parabola or a degenerate conic first as we explain next. The case of a hyperbola follows from $\det(S')$ being continuous on the image of $\tilde\pi$ and  from an ellipse having $\det(S')>0$ and a hyperbola having $\det(S')<0$. For an imaginary ellipse observe that the function $(u_0+u_2)\det(S)$ is also continuous in the coefficients $u_0,\dots,u_5$. For it to change sign, it has to be zero first, in which case either $\det(S)=0$, but then the associated conic is degenerate, or $u_0+u_2=0$. If $u_0=-u_2$, the determinant $\det(S')=-u_0^2-u_1^2/4$ can be negative, which gives a hyperbola (already dealt with above), or it can be zero, in which case $u_0=u_1=u_2=0$, and the corresponding conic is degenerate.

    Let $C_1$ be a circle inside the polygon $P$ and let $C_2$ be a circle containing $Q$.
    These circles exist because $P$ has nonempty interior and $Q$ is bounded.
    By possibly reparametrising the interval in order to consider the subpath $\pi^*:[0,t^*]\to\mathcal C$, the assumptions of~\Cref{prop:limits-of-ellipses} are then satisfied.
    By~\Cref{prop:limits-of-ellipses}, there exists $t_0 \in [0,t^*)$ such that the ellipse $\pi(t_0)$ is not nested between $C_1$ and $C_2$. Hence, the ellipse $\pi(t_0)$ is not nested between $P$ and $Q$. 
\end{proof}

\begin{lemma} \label{lemma:tact_invariant_polygon_containment}
Let $P \subset Q$ be polygons and $E$ an ellipse such that $P \subset \operatorname{conv}(E)$. Let $Q$ be defined by $(x,1)^Tb_j \geq 0$ for $j \in \{1,\ldots,q\}$ and let the defining polynomial for $E$ be $U(x,y)$. Then $T(U,b_j) \leq 0$ for $j \in \{1,\ldots,q\}$ implies that $\operatorname{conv}(E) \subset Q$.
\end{lemma}

\begin{proof}
    Fix $b_j$ and let $H_j$ be the closed half-plane defined by $(x,1)^Tb_j \geq 0$. 
    Then $Q \subset H_j$.
    The condition $T(U,b_j) \leq 0$ asserts that $E$ meets the line $\{x:(x,1)^Tb_j = 0\}$ in at most one point.
    As $E \setminus \{x:(x,1)^Tb_j = 0\}$ is connected, it means that either $E \subset H_j$ or $E \subset H_j^-$ where $H_j^-$ is defined by $(x,1)^Tb_j \leq 0$.
    As $P \subset Q \subset H_j$ and $P \subset \operatorname{conv}(E)$, it follows that $\operatorname{conv}(E) \subset H_j$.
    Since this holds for every edge $b_j$ of $Q$, and $Q = \cap_{j=1}^q H_j$, it follows that $\operatorname{conv}(E) \subset Q$. 
\end{proof}

\begin{proposition} \label{proposition:one_more_incidence}
Let\/ $\mathcal{C}$ be a family of ellipses such that each ellipse in the family contains  $k$ given vertices of\/ $P$ and is tangent to  $\ell$ given edges of\/ $Q$. If\/ $\pi:[0,1] \rightarrow \mathcal{C}$ is a continuous path such that\/ $E_0:=\pi(0)$ is an ellipse nested between $P$ and\/ $Q$ and has no further incidence and tangency relations with $P$ and $Q$, and\/ $E:=\pi(1)$ is an ellipse not nested between $P$ and\/ $Q$, then there exists $t^* \in (0,1)$ such that\/ $E':=\pi(t^*)$ is nested between $P$ and\/ $Q$, and $E'$ has at least one more incidence or tangency with $P$ and\/ $Q$ than $E_0$ does. 
\end{proposition}

\begin{proof}
For each $\pi(t) \in \mathbb{P}^5$, we choose a representative in $\mathbb{R}^6$ such that $u_0 >0$ and the path of representatives is continuous in $\mathbb{R}^6$. Let the defining polynomial of the representative of $\pi(t)$ be $U_t(x,y)$. Let $a_1, \dots , a_p$ be the vertices of $P$  and let $b_1,\ldots,b_q$ define edges of $Q$. Then $P$ is contained in the convex hull of the ellipse $\pi(t)$ if and only if $U_t(a_i) \leq 0$ for $i \in \{1,\ldots,p\}$. If the polygon $Q$ contains the ellipse, then $T(U_t,b_j) \leq 0$ for $j \in \{1,\ldots,q\}$. By~\Cref{lemma:tact_invariant_polygon_containment}, for all $t$ such that $P$ is contained in the convex hull of the ellipse of $\pi(t)$, the opposite implication of the previous sentence also holds.

Let
\[
W=\{t \in [0,1]: \pi(t) \text{ is nested between } P \text{ and } Q\}.
\]
The set $W$ is closed in $\mathbb{R}$ with respect to the Euclidean topology because it is the preimage of $(-\infty,0]$ under the continuous functions $t \mapsto U_t(a_i)$ and $t \mapsto T(U_t,b_j)$. It is also nonempty as $0 \in W$. Every nonempty, closed and bounded set has a maximum element. We take $t^*$ to be this maximum. By construction  $E':=\pi(t^*)$  is nested between $P$ and $Q$. We have $t^*<1$ because $E=\pi(1)$ is not nested between $P$ and $Q$.

We claim that $E'$ has at least one more incidence or tangency with $P$ and $Q$ than $E_0$.  An incidence or tangency with $P$ or $Q$ means that the corresponding inequality $U_t(a_i) \leq 0$ or $T(U_t,b_j) \leq 0$ is satisfied with equality. Let $J \subset \{a_1,\ldots,a_p,b_1,\ldots,b_q\}$ be the index set of inequalities that are satisfied with equality for all ellipses in the family $\mathcal{C}$. Suppose that $E'$ satisfies with equality precisely the inequalities indexed by $J$. This means that for all $t$ in some sufficiently small neighborhood of $t^*$, the strict inequalities persist and the ellipses satisfy with equality precisely the inequalities indexed by $J$, so $\pi(t)$ is nested between $P$ and $Q$ in some neighborhood of $t^*$. This is a contradiction to $t^*$ being the maximum such that $\pi(t^*)$ is nested between $P$ and $Q$. This means that $E'$ satisfies with equality some additional inequalities beyond the index set $J$ and also $t^*>0$ since $E_0$ satisfies with equality precisely the inequalities indexed by $J$.
\end{proof}

Now we are ready to prove the main result in this section:

\begin{proof}[Proof of~\Cref{prop:increasing_incidences}]
Let $E \in \mathcal{C}$ be an ellipse such that $P \subset \operatorname{conv}(E) \subset Q$ and let $C \in \mathcal{C}$ be a conic that is not an ellipse. Then there is a path $\pi:[0,1] \rightarrow \mathcal{C}$ such that $\pi(0)=E$ and $\pi(1)=C$. By~\Cref{cor:path_to_different_conic_implies_non_nestedness}, there is some $t_0\in(0,1)$ such that $\pi(t)$ is an ellipse for $t \in [0,t_0]$ and $\pi(t_0)$ is not nested between $P$ and $Q$.  By~\Cref{proposition:one_more_incidence}, there exists $t^* \in (0,t_0)$ such that $E':=\pi(t^*)$ is nested between $P$ and $Q$ and $E'$ has at least one more incidence or tangency with $P$ and $Q$ than $E$. 
\end{proof}

\section{From infinitely many to finitely many ellipses}\label{sec:infinite_to_finite}

\begin{theorem}[Restatement of Theorem~\ref{thm:main}]
Let $P, Q \subset \mathbb{R}^2$ be polygons and $E \subset \mathbb{R}^2$ an ellipse such that $P \subset \operatorname{conv}(E) \subset Q$. Then there exists an ellipse $E'$ such that $P \subset \operatorname{conv}(E') \subset Q$ and $E'$ has at least five incidence and tangency relations with the polygons $P$ and $Q$, respectively.
\end{theorem}

\begin{proof}
We start by introducing the strategy for the proof before carrying out technicalities for different cases. Let $E$ be an ellipse  such that $P \subset \operatorname{conv}(E) \subset Q$ and $E$ has four incidence and tangency relations with the polygons $P$ and $Q$.  We will consider the family of all conics with given four incidence and tangency relations, and show that this family contains a path-connected component that contains the ellipse $E$ and  a conic that is not an ellipse.   By~ \Cref{prop:increasing_incidences}, there exists an ellipse $E'$ such that $P \subset \operatorname{conv}(E') \subset Q$ and $E'$ has at least five incidence and tangency relations with the polygons~$P$ and~$Q$.

Now suppose $E$ has $k<4$ incidence and tangency relations with $P$ and $Q$. To reduce to the case with four incidence and tangency relations, we fix $4-k$ additional points on $E$, so that we can apply the same argument as for four incidence and tangency relations. This produces a nested ellipse with one additional incidence or tangency with $P$ or $Q$, beyond the original $k$. Repeating this procedure of fixing additional points and reapplying the argument, we increase the number of incidence and tangency relations one at a time until it reaches five.

There are two types of nested polygons $P \subset Q$ we need to consider: the polygons $P$ that are contained in the interior of the polygon $Q$, and the polygons $P$ such that the boundary of $P$ touches the boundary of $Q$. 
In particular, a vertex of $P$ can lie on an edge of $Q$.
The former case corresponds to positive matrices and the latter to matrices with one or more zeros.
At most one vertex of $P$ lies on any given edge of $Q$ and no vertex of $P$ is a vertex of $Q$; otherwise it would not be possible for the ellipse $E$ to exist.

We will first prove the result for the cases with four incidence and tangency relations including at most two line-tangencies with $Q$. The cases with three and four line-tangencies follow from the zero and one line-tangency cases by polarity: By applying an affine transformation, we may assume that the polygons $P$ and $Q$ and the ellipse $E$ contain the origin in their interiors. Hence the polars $P^{\circ}$ and $Q^{\circ}$ are (bounded) polygons and $E^{\circ}$ is an ellipse satisfying $Q^{\circ} \subset \operatorname{conv}(E^{\circ}) \subset P^{\circ}$. A vertex of $P$ lies on $E$ if and only if the corresponding edge of $P^\circ$ is tangent to $E^\circ$, and an edge of $Q$ is tangent to $E$ if and only if the corresponding vertex of $Q^\circ$ lies on $E^\circ$.  If $Q$ has at least three tangencies with $E$, then $P$ has at most one incidence with $E$, which means that $P^{\circ}$ has at most one tangency with $E^{\circ}$, so the statement for these cases follows from the ones for polars. As a point-line incidence becomes a line-point incidence under polarity, a vertex of $P$ on an edge of $Q$ becomes a vertex of $Q^{\circ}$ on an edge of $P^{\circ}$, which is again a zero for the polar pair.

For each configuration of incidence and tangency relations with $P$ and $Q$, we first prove the case when $P$ is contained in the interior of $Q$.
Next, if applicable, we prove the same for the cases when at least one vertex of $P$ lies on an edge of $Q$.
As each incidence between the vertices of $P$ and edges of $Q$ fixes one incidence and one tangency between the ellipse and the polygons $P$ and $Q$
we need to consider only cases where at most two vertices of $P$ lie on edges of $Q$.

For each of the cases with four incidences and tangencies, we show that the path-connected component of $E$ in the family of conics with these four incidences and tangencies contains a conic that is not an ellipse. The first step is to describe the family, or a relevant part of it. 
We do so by either explicitly describing the elements of the path-connected component or via elimination: writing the ideal defining the family as $I$, we eliminate some variables and obtain an elimination ideal $I_l$.
The coordinate projection $\pi_l$ to the remaining coordinates satisfies $\pi_l(\mathbb{V}(I)) \subset \mathbb{V}(I_l)$. 
We then show that the reverse inclusion holds as well by proving that every element in $\mathbb{V}(I_l)$ extends uniquely to an element of $\mathbb{V}(I)$: this proves that $\pi_l(\mathbb{V}(I)) = \mathbb{V}(I_l)$.

In the rest of the proof, a point $P=(P_1,P_2)$ lying on a conic gives the equation 
\begin{equation} \label{eqn:point_incidence}
u_0 P_1^2 + u_1 P_1 P_2 + u_2 P_2^2 + u_3 P_1 + u_4 P_2 + u_5=0.
\end{equation} 
A conic being tangent to $x=0, y=0,$ and $x=1$, respectively, results in the equations:
\begin{equation} \label{eqn:line_tangency}
\begin{split}
&u_4^2-4u_2u_5=0,\quad u_3^2-4u_0u_5=0,\\
&u_1^2-4u_0 u_2-4u_2u_3+2u_1 u_4+u_4^2-4u_2 u_5=0.
\end{split}
\end{equation}

The projective space $\mathbb{P}^5$ contains points that do not correspond to conics, namely when $u_0=u_1=u_2=0$. If $u_3 \neq 0$ or $u_4 \neq 0$, then the resulting geometric object is a line. If only $u_5 \neq 0$, then we get an empty set. We explore when such points of $\mathbb{P}^5$ are solutions to the equations above. A point $u=[u_0:u_1:u_2:u_3:u_4:u_5] \in \mathbb{P}^5$ satisfies the equation~\eqref{eqn:point_incidence} if and only if the point $P$ be belongs to the geometric object corresponding to $u$. Therefore, the only non-conic possibilities are lines that contain $P$. If we have equations for at least three non-collinear points, then there is no line that goes through all the points, and thus also no non-conic solution to the corresponding equations. If we have equations for two points, then the only possibility is the line through these points. 

Next we look which of the lines are solutions to equations coming from the tangency conditions~\eqref{eqn:line_tangency}. By the previous paragraph, we do not have to consider the points with only $u_5 \neq 0$ as our later case analysis will always have some equations for point incidences. We consider the equation $u_4^2-4u_2u_5=0$, which corresponds to tangency to $x=0$. Taking $u_2=0$, we get $u_4=0$, i.e., the corresponding geometric object is a line of the form $x=c$. This is a line parallel to $x=0$. For the other two tangency equations, we similarly get as solutions lines that are parallel to the tangency line. From this analysis it follows that in the case of a total of four point incidence and line tangency conditions, out of which at least two are point incidences, the only setting where we can get non-conic solutions to the corresponding equations are when we have two parallel lines and two points such that the line through them is parallel to the two lines as well. In all other such cases, there are no non-conic solutions.

 \textbf{Four points, no lines, no zeros:} Let the incidences be given by points $A,B,C,D$. Since $A,B,C,D$ lie on an ellipse, no three of them are collinear, so the four linear equations they impose on the coefficients $u_0,\dots,u_5$ of $\mathbb{P}^5$ are independent. The solution set of this system is therefore a one-dimensional linear space in $\mathbb{P}^5$; in particular, it is path-connected.
This linear space also contains the pair of lines $l_{AB}$ through points $A,B$ and $l_{CD}$ through $C,D$, which is a degenerate conic.

\textbf{Three points, one line, no zeros:} Let the incidences be given by three points $A,B,C$ and the tangency by one line $d$ such that none of $A,B,C$ lies on $d$. We can assume that after an affine transformation the line $d$ is given by $x=0$. The ellipse $E$ remains an ellipse under the affine transformation. The conics satisfying these incidences and tangencies form the solution set of
\begin{equation} \label{eqn:family_of_incident_conics_3_1}
\begin{cases}
u_4^2-4u_2u_5=0,\\
u_0 A_1^2 + u_1 A_1 A_2 + u_2 A_2^2 + u_3 A_1 + u_4 A_2 + u_5=0,\\
u_0 B_1^2 + u_1 B_1 B_2 + u_2 B_2^2 + u_3 B_1 + u_4 B_2 + u_5=0,\\
u_0 C_1^2 + u_1 C_1 C_2 + u_2 C_2^2 + u_3 C_1 + u_4 C_2 + u_5=0.
\end{cases}
\end{equation}
For each $(u_2,u_4,u_5)$ satisfying $u_4^2-4u_2u_5=0$, the three linear equations determine $(u_0,u_1,u_3)$ uniquely, because their coefficient matrix
\[
\begin{pmatrix}
A_1^2 & A_1 A_2 & A_1\\
B_1^2 & B_1 B_2 & B_1\\
C_1^2 & C_1 C_2 & C_1
\end{pmatrix}
\]
is nonsingular. Indeed, its determinant is
\[
\det
\begin{pmatrix}
A_1^2 & A_1 A_2 & A_1\\ 
B_1^2 & B_1 B_2 & B_1\\ 
C_1^2 & C_1 C_2 & C_1
\end{pmatrix}
= A_1 B_1 C_1 (-A_2 B_1 + A_1 B_2 + A_2 C_1 -B_2 C_1 - A_1 C_2 + B_1 C_2),
\]
which is nonzero: $A_1,B_1,C_1$ are nonzero because $E$ meets $x=0$ only at the tangency point on $x=0$ and $A,B,C$ are point-incidences with the ellipse $E$ distinct from it, and the second factor is nonzero because $A,B,C$ lie on the ellipse $E$ and hence cannot be collinear.

Since $\det\begin{pmatrix}0&0&-2\\0&1&0\\-2&0&0\end{pmatrix}=-4\neq0$ and the conic $u_4^2-4u_2u_5=0$ has the real point $[u_2:u_4:u_5]=[1:0:0]$, it is a nonsingular real conic in $\mathbb{P}^2$, hence topologically a circle, in particular path-connected. The extension $(u_2,u_4,u_5)\mapsto(u_0,u_1,u_2,u_3,u_4,u_5)$ above is linear in $u_2,u_4,u_5$, and its inverse is the coordinate projection onto $(u_2,u_4,u_5)$. Hence the solution set~\eqref{eqn:family_of_incident_conics_3_1} in $\mathbb{P}^5$ is homeomorphic to $u_4^2-4u_2u_5=0$ in $\mathbb{P}^2$, and therefore path-connected as well.

This space also contains a degenerate conic. Since $A,B,C$ are not collinear, at most one of the lines $l_{AB},l_{AC},l_{BC}$ is parallel to $d$. Relabel $A,B,C$ if necessary so that $l_{AB}$, the line through $A$ and $B$, is not parallel to $d$, and set $P=l_{AB}\cap d$. Let $l_{CP}$ be the line through $C$ and $P$. Writing $l_{AB}: a_1x+b_1y+c_1=0$ and $l_{CP}: a_2x+b_2y+c_2=0$, the conic $l_{AB}\cup l_{CP}$ has $u_4^2-4u_2u_5=(b_1c_2-b_2c_1)^2$, which vanishes because $l_{AB}$ and $l_{CP}$ meet $d$ at the common point $P$. Hence $l_{AB}\cup l_{CP}$ is a degenerate conic satisfying~\eqref{eqn:family_of_incident_conics_3_1}.

\textbf{Three points, one line, one zero:} Let the incidences be given by three points $A,B,C$ and the tangency by a line $d$ such that the point $A$ lies on the line $d$.
    After applying an affine transformation, we may assume that the line $d$ is $x=0$.
    Moreover, we may assume that $A=(0,0)$ after applying a further translation in the $y$-direction.
    The ellipse remains an ellipse after this affine transformation.
    Since $A$ is the only point of the ellipse on the line $x=0$, it follows that $B_1,C_1 \neq 0$.

First, we assume that $B_1-C_1 \neq 0$.
The defining equation of the ellipse $E$ satisfies $u_2 \neq 0$ by ~\Cref{rmk:real_ellipse_nonzero_coefficients}.
    The degenerate conic we construct below also satisfies $u_2 \neq 0$, so we will restrict to the affine chart $u_2 \neq 0$.
We get the following system of equations:
    \begin{equation*}
        \begin{cases}
            u_4^2 - 4 u_2 u_5 =0, \\
            u_5 = 0, \\
            B_1^2 u_0 + B_1 B_2 u_1 + B_2^2 u_2 + B_1 u_3 + B_2 u_4 + u_5 = 0, \\
            C_1^2 u_0 + C_1 C_2 u_1 + C_2^2 u_2 + C_1 u_3 + C_2 u_4 + u_5 = 0, \\
            u_2 - 1 = 0.
        \end{cases}
    \end{equation*}
    Denote by $V_{A,B,C,d}$ the solution set of this system
    and let $(u_0,\dots, u_5) \in V_{A,B,C,d}$.
    The second and fifth equations immediately imply that $u_2=1$ and $u_5 = 0$.
    Substituting these into the first equation gives $u_4 = 0$.
    Substituting all of these in turn to the third and fourth equations gives the equations $B_1^2 u_0 + B_1 B_2 u_1 + B_1 u_3 + B_2^2 = 0$ and $C_1^2 u_0 + C_1 C_2 u_1 + C_1 u_3 + C_2^2  = 0.$
    They are linear in $(u_0, u_3)$ with coefficient matrix $\begin{pmatrix}
        B_1^2 & B_1 \\
        C_1^2 & C_1
    \end{pmatrix}$ whose determinant $B_1 C_1 (B_1-C_1)$ is nonzero, as $B_1,C_1 \neq 0$ and by assumption $B_1 -C_1 \neq 0$. Hence, $u_0$ and $u_3$ can be expressed as a linear function of $u_1$.
    Combining these observations, every element of $V_{A,B,C,d}$ is of the form $(u_0(u_1), u_1, 1, u_3(u_1), 0, 0)$, where $u_0(u_1)$ and $u_3(u_1)$ are linear functions, for $u_1 \in \R$.
    Conversely, every element of this form satisfies all defining equations of the family $V_{A,B,C,d}$.
    Hence, the family $V_{A,B,C,d}$ is the line $\{(u_0(u_1), u_1, 1, u_3(u_1), 0, 0) \mid u_1 \in \R\}$, that is, path-connected.

    The family $V_{A,B,C,d}$ also contains a pair of lines.
        Let $l_{AB}$ be the line through the points $A$ and $B$ given by $B_2x-B_1y= 0$ and let $l_{AC}$ be the line through $A$ and $C$ given by $C_2x-C_1y=0$.
    Then the pair of lines $l_{AB} \cup l_{AC}$ has $u_2 = B_1 C_1 \neq 0$; rescaling so that $u_2 = 1$ places it in $V_{A,B,C,d}$.
    Since $u_2 \neq 0$ for the original ellipse $E$ as well, it follows that $E \in V_{A,B,C,d}$, whose single path-connected component contains both the ellipse $E$ and this non-ellipse conic.

    Now, assume $B_1 -C_1 = 0$ and write $C=(B_1,C_2)$.
    We consider the same system without restricting to the affine chart $u_2 \neq 0$ and substitute $C_1 = B_1$.
    The conics tangent to $x = 0$ passing through the points $A, B$, and $C$ are solutions to the following system:
    \begin{equation*}
        \begin{cases}
            u_4^2 - 4u_2 u_5 =0, \\
            u_5 = 0, \\
            B_1^2 u_0 + B_1 B_2 u_1 + B_2^2 u_2 + B_1 u_3 + B_2 u_4 + u_5 = 0, \\
            B_1^2 u_0 + B_1 C_2 u_1 + C_2^2 u_2 + B_1 u_3 + C_2 u_4 + u_5 = 0.
        \end{cases}
    \end{equation*}
    Let $V_{A,B,c,d}$ denote the solution set of this system in $\mathbb{R}^6$ and let $(u_0, \dots, u_5)$ be an element of $V_{A,B,c,d}$.
    The second equation implies that $u_5 = 0$; substituting this to the first equation gives $u_4 ^2 = 0$, that is, $u_4 =0$.
    With these, the equations three and four reduce to 
    linear equations in the variables $u_0,u_1,u_2,u_3$.
    Subtracting the fourth equation from the third and factoring yields $(B_2 - C_2) (B_1 u_1 + (B_2 + C_2) u_2) = 0$, which, since $B_2 \neq C_2$, is equivalent to $B_1 u_1 + (B_2 + C_2) u_2 = 0$.
    Since $B_1 \neq 0$ it follows that $u_1$ depends linearly on the value of the variable $u_2$.
    Substituting the expression we get for $u_1$ to the third equation yields an equation linear on $u_0,u_2,u_3$.
    Hence, using the fact that $B_1 \neq 0$,  $u_0$ can be expressed uniquely as a linear function of $u_2,u_3$.
    Thus, $(u_0,\dots, u_5) = (u_0(u_2,u_3), u_1(u_2), u_2,u_3,0,0)$.
    Conversely, any point of the above form satisfies the defining equations of $V_{A,B,c,d}$.
    Hence, $V_{A,B,c,d} = \{(u_0(u_2,u_3), u_1(u_2), u_2,u_3,0,0) \mid u_2, u_3 \in \R\}$ is a 2-dimensional linear subspace in $\mathbb{R}^6$, and the corresponding family of conics is the projective line $\mathbb{P}(V_{A,B,c,d})\subset\mathbb{P}^5$.
    In particular, it is path-connected.

    The family $\mathbb{P}(V_{A,B,c,d})$ contains the degenerate conic consisting of a pair of lines. 
    Let $l_{AB}$ be the line through the points $A$ and $B$ given by $B_2x-B_1y= 0$ and let $l_{AC}$ be the line through $A$ and $C$ given by $C_2x-B_1y=0$.
    The pair of lines $l_{AB} \cup l_{AC}$ is a degenerate conic contained in $\mathbb{P}(V_{A,B,c,d})$.
    Hence, $\mathbb{P}(V_{A,B,c,d})$ contains both the ellipse and the degenerate conic.

\textbf{Two points, two intersecting lines, no zeros:} Let the incidences be given by two points $A,B$ and the tangencies by two intersecting lines $c,d$ such that $A,B$ do not lie on $c,d$. After applying an affine transformation, we can assume that the line $c$ is given by $x=0$ and the line $d$ is given by $y=0$. Since $A,B$ lie on the ellipse $E$ that is tangent to the lines $x=0$ and $y=0$, none of the coordinates of $A,B$ are zero.

Since the ellipse $E$ is tangent to $x=0$ and $y=0$, it cannot contain the origin, so $u_5\neq0$ for $E$. The non-ellipse conic constructed below will also lie in this chart, so we will work in the affine chart $u_5\neq0$. The conics with the four incidences and tangencies in this chart are solutions to the following system of equations:
\begin{equation*} 
\begin{cases}
u_4^2-4u_2u_5=0,\\
u_3^2-4u_0u_5=0,\\
u_5-1=0,\\
u_0 A_1^2 + u_1 A_1 A_2 + u_2 A_2^2 + u_3 A_1 + u_4 A_2 + u_5=0,\\
u_0 B_1^2 + u_1 B_1 B_2 + u_2 B_2^2 + u_3 B_1 + u_4 B_2 + u_5=0.
\end{cases}
\end{equation*}
We denote by $V_{A,B,c,d}$ the solution set of this system.
Eliminating the variables $u_0,u_1,u_2,u_5$ from the system gives the following polynomial, whose zero set in the $(u_3,u_4)$-plane we denote by $\tilde{C}$:
\begin{align*}
& u_3^2 A_1 A_2 B_1^2-u_3^2 A_1^2 B_1 B_2-u_4^2 A_2^2 B_1 B_2+u_4^2 A_1 A_2 B_2^2+\\
&4 u_3 A_1 A_2 B_1+4 u_4 A_1 A_2 B_2-4 u_3 A_1 B_1 B_2-4 u_4 A_2 B_1 B_2+4 A_1 A_2-4 B_1 B_2 =0.
\end{align*}
This elimination can be carried out by the same procedure regardless of $A,B$: after substituting $u_5=1$, the first equation gives $u_2=u_4^2/4$ and the second equation gives $u_0 = u_3^2/4$. The coefficient of $u_1$ in each of the last two equations, $A_1A_2$ and $B_1B_2$ respectively, is nonzero because no coordinate of $A,B$ is zero, so substituting $u_0,u_2,u_5$ and solving each equation for $u_1$ gives two expressions for $u_1$ in terms of $u_3,u_4$; equating them, and clearing the denominators $A_1A_2$ and $B_1B_2$, eliminates $u_1$ and yields the single polynomial relation above.

The matrix associated to this quadratic polynomial is
\[
\tilde{Q}=
\begin{pmatrix}
A_1 A_2 B_1^2 - A_1^2 B_1 B_2  &  0  &  2 A_1 A_2 B_1 - 2 A_1 B_1 B_2\\   
0  &  A_1 A_2 B_2^2 - A_2^2 B_1 B_2  &  2 A_1 A_2 B_2 - 2 A_2 B_1 B_2\\   
2 A_1 A_2 B_1 - 2 A_1 B_1 B_2  & 
2 A_1 A_2 B_2 - 2 A_2 B_1 B_2  & 
4 A_1 A_2 - 4 B_1 B_2
\end{pmatrix}.
\]
One can compute that
\[
\det(\tilde{Q})=0 \text{ and } \det(\tilde{Q}_{12,12})=-A_1 A_2 B_1 B_2 (A_2 B_1 - A_1 B_2)^2.
\]
There are two possibilities for $\tilde{C}$, according to whether $A_2 B_1 - A_1 B_2$ vanishes. First assume $A_2 B_1 - A_1 B_2 \neq 0$. Since $E$ is tangent to $x=0$ and $y=0$, it lies entirely in one orthant, so $A,B\in E$ lie in the same orthant; hence $A_1 A_2 B_1 B_2 >0$. Since $(A_2B_1-A_1B_2)^2>0$, it follows that $\det(\tilde{Q}_{12,12}) <0$, so $\tilde{C}$ is two real intersecting lines. Now assume $A_2 B_1 - A_1 B_2=0$. Then the coefficients of $u_3^2,u_3u_4,u_4^2$ are zero. Since the coefficient of $u_3$, which is $4A_1 B_1(A_2-B_2)$, and the coefficient of $u_4$, which is $4A_2(A_1-B_1)B_2$, cannot both vanish ($A,B$ are distinct), then $\tilde{C}$ is a line in the $(u_3,u_4)$-plane. In either case, $\tilde{C}$ is path-connected: the two intersecting lines meet at a point, and the single line is trivially connected.

Since the extension $(u_3,u_4) \mapsto (u_0,u_1,u_2,u_3,u_4,u_5)$ is a polynomial function in $u_3,u_4$, with inverse the coordinate projection onto $(u_3,u_4)$, this extension is a homeomorphism between $\tilde{C}$ and $V_{A,B,c,d}$. Since $\tilde{C}$ is path-connected,  it follows that $V_{A,B,c,d}$ is path-connected.

Since every line contains a point with $u_3=0$ or with $u_4=0$, $\tilde{C}$ contains such a point. If $u_4=0$, then $u_4^2-4u_2u_5=0$ and $u_5=1$ imply $u_2=0$. The matrix of the corresponding conic in $V_{A,B,c,d}$ is
\[
Q=
\begin{pmatrix}
u_0 & u_1/2 & u_3/2\\
u_1/2 & 0 & 0\\
u_3/2 & 0 & 1
\end{pmatrix}.
\]
Then $\det(Q) = -u_1^2/4$ and $\det(Q_{12,12})=-u_1^2/4$. If $\det(Q)=0$, then it is a degenerate conic, and if $\det(Q) \neq 0$, then it is a hyperbola. The case $u_3=0$ gives the same $\det(Q)$ and $\det(Q_{12,12})$, and hence leads to the same conclusion. In either case, this gives a conic in $V_{A,B,c,d}$ that is not an ellipse. Since $V_{A,B,c,d}$ has a single path-connected component, this component contains both $E$ and this non-ellipse conic.

\textbf{Two points, two intersecting lines, one zero:} Let the incidences be given by points $A,B$ and the tangencies by two intersecting lines $c,d$ such that the point $A$ lies on line $c$.
    After applying an affine transformation, we may assume that the line $c$ is given by $x= 0$ and the line $d$ is given by $y=0$. 
    Moreover, we may assume that the point $A=(0,1)$ by applying a further scaling in the $y$-direction.
    Since the point $B$ does not lie on the lines $c$ or $d$ it follows that $B=(B_1,B_2)$ where $B_1,B_2 \neq 0$.
    The ellipse stays an ellipse after an affine transformation.

    We may restrict to the affine chart $u _5 \neq 0 $ since no nondegenerate conic tangent to the lines $x=0$ and $y=0$ can pass through the origin.
    The non-ellipse conic we construct is also in this chart, so we need not leave it.
    The conics tangent to the lines $x = 0$ and $y = 0$ passing through the points $A$ and $B$ are the solutions to the following system of equations:
    \begin{equation*}
        \begin{cases}
            u_4^2-4 u_2 u_5 = 0,\\
            u_3^2 - 4 u_0 u_5 = 0,\\
            u_2 + u_4 +u_5 = 0, \\
            B_1 ^2 u_0 + B_1 B_2 u_1 + B_2^2 u_2 + B_1 u_3 + B_2 u_4 +u_5 = 0, \\
            u_5 - 1 =0.
        \end{cases}
    \end{equation*}
    Let $V_{A,B,c,d}$ denote the solution set of this system
    and let $(u_0,\dots, u_5) \in V_{A,B,c,d}$.
    After substituting $u_5 = 1$ to the second equation, the variable $u_0$ is uniquely determined from the variable $u_3$ by $u_0 = u_3^2/4$.
    After substituting $u_5 = 1$ to the third equation $u_2 +u_4 + u_5 =0$, the variable $u_2$ is uniquely determined as a linear function of $u_4$ from $u_2 = -u_4-1$.
    Substituting these expressions for $u_2$ and $u_5$ to the first equation gives $(u_4+2)^2 =0$, that is, $u_4 = -2$.
    The fourth equation has constant coefficient $B_1 B_2 \neq 0$ for the variable $u_1$, so, together with the expressions for $u_0$, $u_2$, $u_4$ and $u_5$, it uniquely determines $u_1$ as a linear function of $u_3$.
    Combining these observations, every element of $V_{A,B,c,d}$ is of the form $(u_0(u_3), u_1(u_3), 1,u_3, -2, 1)$ for $u_3 \in \R$.
    Conversely, every element of this form satisfies all defining equations of the family $V_{A,B,c,d}$.
    Hence, $V_{A,B,c,d}=\{(u_0(u_3), u_1(u_3), 1,u_3, -2, 1) \mid u_3 \in \R\}$, which is homeomorphic to $\mathbb{R}$, and thus path-connected.
    
    Lastly, we show that the family $V_{A,B,c,d}$ also contains a degenerate conic. Consider the double line through $A$ and $B$. Since the line $l_{AB}$ through $A$ and $B$ does not go through the origin, the equation defining the corresponding double line has nonzero constant coefficient and we can rescale it to have $u_5=1$. The double line clearly satisfies the equations for the point-incidences. One can easily check that every double line also satisfies the tact invariants. Hence the double line is in the family $V_{A,B,c,d}$.
    Thus the path-connected family $V_{A,B,c,d}$ contains both the ellipse $E$ and the non-ellipse conic.

\textbf{Two points, two intersecting lines, two zeros:} Let the incidences be given by two points $A,B$ and the tangencies by two intersecting lines $c,d$ such that $A$ lies on $c$ and $B$ lies on $d$.
    After applying an affine transformation, we may assume that the line $c$ is given by $x = 0$, the line $d$ is given by $y=0$, and the points are $A=(0,1)$ and $B=(1,0)$.
    The ellipse $E$ remains an ellipse after an affine transformation.
    
    Since the ellipse $E$ is tangent to the lines $x=0$ and $y=0$, it cannot contain the origin, so $u_5 \neq 0$ for $E$.
    The non-ellipse conic constructed below will also have $u_5 \neq 0$, so we will work in the affine chart $u_5 \neq 0$.
    The conics with the four incidences and tangencies with $A,B,c,d$ in this chart are solutions to the following system of equations:
    \begin{equation*}
        \begin{cases}
            u_4^2-4 u_2u_5 =0,\\
            u_3^2-4 u_0u_5=0,\\
            u_5 -1=0,\\
            u_0+ u_3  + u_5= 0,\\
            u_2 + u_4 + u_5 = 0.
        \end{cases}
    \end{equation*}
    Let $V_{A,B,c,d}$ denote the solution set of the system of these equations 
    and let $(u_0,\dots,u_5) \in V_{A,B,c,d}$.
    After substituting $u_5 = 1$, the first equation gives $u_2 = u_4^2/4$ and the second equation gives $u_0 = u_3^2/4$.
    Substituting these expressions for $u_0,u_2$ and $u_5$ into the fourth and fifth equations and clearing out denominators leaves the equations $(u_3+2) ^2 =0$ and $(u_4+2)^2 =0$, that is, $u_3=-2$ and $u_4 = -2$.
    Therefore, every element of $V_{A,B,c,d}$ is of the form $(1,u_1,1,-2,-2,1)$ for $u_1 \in \R$.
    Conversely, every element of this form satisfies the defining equations of $V_{A,B,c,d}$.
    Hence, $V_{A,B,c,d}$ is the line $\{(1,u_1,1,-2,-2,1) \mid u_1 \in \R\}$ which is path-connected.

    The double line $(x+y-1)^2=0$ through the points $A$ and $B$ is a degenerate conic that belongs to $V_{A,B,c,d}$.
    Since $V_{A,B,c,d}$ is path-connected, it contains both the ellipse $E$ and the degenerate conic.

\textbf{Two points, two parallel lines, no zeros:}  Let the incidences be given by two points $A,B$ and the tangencies by two parallel lines $c,d$ such that $A,B$ do not lie on $c,d$. After applying an affine transformation, we can assume that the line $c$ is given by $x=0$ and the line $d$ is given by $x=1$. Moreover, we can assume that $A$ is of the form $(A_1,0)$ by applying a further translation in the $y$-direction. Since $A,B$ are contained in the ellipse $E$ that is tangent to the lines $x=0$ and $x=1$,  the first coordinates of $A,B$ satisfy $0<A_1,B_1<1$. Moreover, for the ellipse $E$ we have $u_2 \neq 0$ by~\Cref{rmk:real_ellipse_nonzero_coefficients}.

First we consider the case $A_1 - B_1 \neq 0$ and we restrict to the affine chart $u_2 \neq 0$. The conics tangent to $x=0$ and $x=1$ through the points $A,B$ satisfy:
\begin{equation} \label{eqn:family_of_incident_conics_2pts_2par}
\begin{cases}
u_4^2-4u_2u_5=0,\\
u_1^2-4u_0 u_2-4u_2u_3+2u_1 u_4+u_4^2-4u_2 u_5=0,\\
u_2-1=0,\\
u_0 A_1^2 + u_3 A_1 + u_5=0,\\
u_0 B_1^2 + u_1 B_1 B_2 + u_2 B_2^2 + u_3 B_1 + u_4 B_2 + u_5=0.
\end{cases}
\end{equation}
Specializing to $u_2=1$ and eliminating the variables $u_0,u_3,u_5$ from the system gives the following polynomial, whose zero set in the $(u_1,u_4)$-plane we denote by $\tilde{C}$:
\begin{align*}
&u_1^2 A_1^2 B_1 - u_1^2 A_1 B_1^2 + 2u_1u_4 A_1^2 B_1 - 2u_1u_4 A_1 B_1^2 + u_4^2 A_1^2 - u_4^2 A_1 - u_4^2 B_1^2 + u_4^2 B_1\\
&+4u_1 A_1^2 B_1 B_2 - 4u_1 A_1 B_1 B_2 + 4u_4 A_1^2 B_2 - 4u_4 A_1 B_2 + 4A_1^2 B_2^2 - 4A_1 B_2^2=0.
\end{align*}
This elimination can be carried out by the same procedure regardless of $A,B$: substituting $u_2=1$ into the first equation gives $u_5=u_4^2/4$. Substituting both of these in turn into the point equations for $A,B$ leaves two equations linear in $(u_0,u_3)$, with coefficient matrix $\begin{pmatrix}A_1^2 & A_1\\ B_1^2 & B_1\end{pmatrix}$, whose determinant $A_1B_1(A_1-B_1)$ is nonzero because $A_1,B_1\neq0$ and $A_1\neq B_1$. Hence $(u_0,u_3)$ are uniquely determined as linear functions of $u_1,u_4$. Substituting $u_2=1$, $u_5=u_4^2/4$, and these expressions for $u_0,u_3$ into the second equation  then yields the polynomial relation in $u_1,u_4$ above. 

The matrix associated to this quadratic polynomial is
\[
\tilde{Q}=\begin{pmatrix}
A_1 B_1 (A_1 - B_1) & A_1 B_1 (A_1 - B_1) & 2 A_1 B_1 B_2 (A_1 - 1) \\
A_1 B_1 (A_1 - B_1) & (A_1 - B_1)(A_1 + B_1 - 1) & 2 A_1 B_2 (A_1 - 1) \\
2 A_1 B_1 B_2 (A_1 - 1) & 2 A_1 B_2 (A_1 - 1) & 4 A_1 B_2^2 (A_1 - 1)
\end{pmatrix}.
\]
One can compute that
\[
\det(\tilde{Q})=0 \text{ and }\det(\tilde{Q}_{12,12})=-A_1 B_1 (1 - A_1)(1 - B_1)\, (A_1 - B_1)^2.
\]
Since $A_1,B_1$ are strictly between 0 and 1, $-A_1 B_1 (1 - A_1)(1 - B_1)<0$. By assumption $A_1 - B_1 \neq 0$, so $\det(\tilde{Q}_{12,12})<0$, and we get two real intersecting lines.

Since $(u_0,u_1,u_2,u_3,u_4,u_5)$ is a polynomial, in particular continuous, function of $(u_1,u_4)$, with inverse the coordinate projection onto $(u_1,u_4)$, this extension is a homeomorphism between $\tilde{C}$ and the solution set of \eqref{eqn:family_of_incident_conics_2pts_2par}. Since $\tilde{C}$ consists of two lines meeting at a point, it is path-connected, and hence so is the solution set of \eqref{eqn:family_of_incident_conics_2pts_2par}.

This family contains the degenerate conic given by the double line through $A$ and $B$: since $A_1\ne B_1$, the line $l_{AB}$ is not vertical, so writing $l_{AB}: ax+by+c=0$ we have $b\ne0$, and the double line $(ax+by+c)^2=0$ has $u_2=b^2\ne0$; rescaling so that $u_2=1$ places it in the family above.  This conic is degenerate, hence not an ellipse. Since $u_2\neq0$ for $E$ as well, $E$ lies in the same family, whose single path-connected component therefore contains both $E$ and this non-ellipse conic.

Now assume $A_1 - B_1=0$. Then $B_2 \neq 0$, because $A$ and $B$ are different. We consider the same system without restricting to the affine chart $u_2 \neq 0 $, substituting $B_1$ by $A_1$. Decomposing the corresponding ideal over $\mathbb{Q}$, where we consider $A_1,B_2$ as ring variables, gives three components
\begin{align*}
&I_1=\langle B_2,4u_5u_2-u_4^2,4u_3u_2+4u_0u_2-u_1^2-2u_1u_4,u_3A_1+u_0A_1^2+u_5\rangle,\\
& I_2=\langle u_2,u_4,u_1,u_3A_1+u_0A_1^2+u_5\rangle,\\
& I_3=\langle u_1A_1+u_4+u_2B_2,4u_5u_2-u_4^2,4u_3u_2+4u_0u_2-u_1^2-2u_1u_4,u_3A_1+u_0A_1^2+u_5,\\
  &4u_0A_1^3-4u_0A_1^2-4u_5A_1+4u_5-2u_4A_1B_2+2u_4B_2+u_2B_2^2,\\
  &4u_0u_4A_1^2-4u_0u_4A_1-4
       u_5u_1-4u_5u_4-u_1u_4B_2-2u_4^2B_2,\\
       &16u_3u_5-16u_0u_5A_1^2+32u_0u_5A_1-4u_0u_2A_1B_2^2+4u_0u_2B_2^2+16u_5^2+8u_5u_1B_2+8u_5u_4B_2-u_1^2B_2^2-2u_1u_4B_2^2\rangle
\end{align*}
with associated varieties $V_1,V_2,V_3$. Although $I_1,I_2,I_3$ might not be the minimal primes of $I$ for all parameter values, then their varieties cover the variety of $I$.
The first component $I_1$ is irrelevant for us since $B_2 \neq 0$. The second component $I_2$ defines a one-dimensional linear space in $\mathbb{P}^5$. This component also contains a non-conic solution $[0:0:0:1:0:-A_1]$. The conic solutions make up $\mathbb{P}^1$ minus a point, which is path-connected. To understand the third component $I_3$, we eliminate the variables $u_0,u_2,u_3$. This elimination can be carried out by the same procedure regardless of the parameters $A_1,B_2$: the first generator $u_1A_1+u_4+u_2B_2$ has coefficient $B_2\ne0$ for $u_2$, so it determines $u_2$ as a linear function of $u_1,u_4$. The fifth generator has coefficient $4A_1^2(A_1-1)\ne0$ for $u_0$, since $A_1\notin\{0,1\}$, so together with the expression for $u_2$ it determines $u_0$. The fourth generator has coefficient $A_1\ne0$ for $u_3$, so together with the expression for $u_0$ it determines $u_3$. Substituting these expressions for $u_0,u_2,u_3$ into the remaining four generators, each reduces to a scalar multiple of a single quadratic expression in $u_1,u_4,u_5$. The scalars depend only on $A_1,B_2$, for two generators they vanish identically at the special value $A_1=1/2$, and they are nonzero for all other $A_1\notin\{0,1\}$ and $B_2\ne0$. Therefore, for every choice of parameters with $A_1\notin\{0,1\}$ and $B_2\ne0$, the elimination ideal is generated by this one quadric. The result is
\[
\langle 4u_5u_1A_1+4u_5u_4+u_4^2B_2 \rangle,
\]
which is a quadratic polynomial with the matrix
\[
\tilde{Q} = \begin{pmatrix} 0 & 0 & 2A_1 \\ 0 & B_2 & 2 \\ 2A_1 & 2 & 0 \end{pmatrix}.
\]

Since $\det(\tilde{Q})=-4A_1^2B_2\neq0$ and the conic $4u_5u_1A_1+4u_5u_4+u_4^2B_2=0$ has the real point $[u_1:u_4:u_5]=[0:0:1]$, it is a nonsingular real conic in $\mathbb{P}^2$, hence topologically a circle, in particular path-connected. Since the extension $(u_1,u_4,u_5)\mapsto(u_0,u_1,u_2,u_3,u_4,u_5)$ above is linear, with inverse the coordinate projection onto $(u_1,u_4,u_5)$, it is a homeomorphism onto its image. By substituting the possible non-conic solution $[0:0:0:1:0:-A_1]$ to the generators of $I_3$, one can check that it does not belong to this component.  Hence the corresponding family of conics is path-connected as well.

Moreover, both the second and third components contain the point corresponding to the double line through $A$ and $B$, given by $(x-A_1)^2=x^2-2A_1x+A_1^2$. Since this point lies in both components, and each component is path-connected, their union is path-connected. Hence this path-connected union contains both $E$ and this degenerate, non-ellipse conic.

    \textbf{Two points, two parallel lines, one zero:} Let the incidences be given by the points $A,B$ and the tangencies by two parallel lines $c, d$ such that $A$ lies on $c$.
    After applying an affine transformation, we may assume that the line $c$ is given by $x=0$ and the line $d$ is given by $x=1$.
    After applying a further translation in the $y$-direction, we may assume that $A=(0,0)$.
    Then $B=(B_1,B_2)$ where $B_1 \neq 0$ and $B_1 \neq 1$.
    Moreover, since $B$ is contained in the ellipse $E$ that is tangent to the lines $x= 0$ and $x=1$, the first coordinate of $B$ satisfies $0<B_1<1$.
    The ellipse $E$ remains an ellipse after the affine transformation.

    For the original ellipse $E$ it must hold that $u_2 \neq 0$ by~\Cref{rmk:real_ellipse_nonzero_coefficients}.
    The degenerate conic we construct will also be contained in this chart. Hence, we work in the affine chart $u_2 \neq 0$.
    The conics tangent to the lines $x=0$ and $x= 1$ passing through the points $A$ and $B$ satisfy the following system of equations:
    \begin{equation*}
        \begin{cases}
            u_4^2 - 4 u_2 u_5 = 0, \\
            u_1^2 - 4 u_0 u_2 - 4 u_2 u_3 + 2 u_1 u_4 + u_4^2 - 4 u_2 u_5=0, \\
            u_5=0, \\
            B_1^2 u_0 + B_1 B_2 u_1 + B_2^2 u_2 + B_1 u_3 + B_2 u_4 + u_5 = 0,\\
            u_2-1=0.
        \end{cases}
    \end{equation*}
    Let $V_{A,B,c,d}$ denote the solution set to this system of equations
    and let $(u_0,\dots,u_5) \in V_{A,B,c,d}$.   The line through $A$ and $B$ cannot be parallel to $x=0$ and hence there are no non-conic solutions.
    From the third and fifth equations it follows that $u_5=0$ and $u_2 = 1$.
    Substituting these into the first equation gives $u_4 = 0$.
    Substituting $u_2,u_4$ and $u_5$ in turn into the second and fourth equations leaves two equations linear in $u_0,u_3$ and quadratic in $u_1$. The coefficient matrix of $(u_0,u_3)$ is $\begin{pmatrix}
            -4 & -4 \\
            B_1^2 & B_1
        \end{pmatrix}$, which is nonsingular whenever $B_1 \neq 0,1$.
    In particular, the variables $u_0$ and $u_3$ are uniquely determined as quadratic functions of the variable $u_1$.
    Combining these observations, every element of $V_{A,B,c,d}$ is of the form $(u_0(u_1), u_1, 1, u_3(u_1), 0,0)$.
    Conversely any point of this form satisfies the defining equations of the family $V_{A,B,c,d}$.
    Hence, $V_{A,B,c,d}=\{(u_0(u_1), u_1, 1, u_3(u_1), 0,0) \mid u_1 \in \R \}$, which is homeomorphic to $\mathbb{R}$, and thus path-connected.

    The family $V_{A,B,c,d}$ contains the degenerate conic given by the double line through the points $A$ and $B$.
    Since $B_1 \neq 0$, the line $l_{AB}$ through the points $A$ and $B$ is not vertical and, since it passes through the origin, we have $l_{AB} \colon ax+by=0$ where $b \neq 0$, and the double line $(ax+by)^2$ has $u_2 = b^2 \neq 0$; rescaling so that $u_2 = 1$ places it in $V_{A,B,c,d}$.
    This conic is degenerate, hence, not an ellipse.
    Since $u_2 \neq 0$ for $E$ as well, $E$ lies in $V_{A,B,c,d}$, whose single path-connected component contains both $E$ and this non-ellipse conic.

    \textbf{Two points, two parallel lines and two zeros:} Let the incidences be given by two points $A,B$ and the tangencies by two parallel lines $c,d$ such that $A$ lies on the line $c$ and $B$ lies on the line $d$.
    After applying an affine transformation we may assume that the line $c$ is given by $x=0$ and the line $d$ is given by $x=1$.
    Moreover, we can assume that $A=(0,0)$ by applying a further translation in the $y$-direction.
    The point $B$ is of the form $B=(1,B_2)$ for $B_2 \in \R$.
    The original ellipse $E$ remains an ellipse after an affine transformation.
    Lastly, for the original ellipse $E$ the term $u_2 \neq 0$ by~\Cref{rmk:real_ellipse_nonzero_coefficients}.
    Since the degenerate conic we construct below is contained in this chart as well, we will work in the affine chart $u_2 \neq 0 $.

    The conics tangent to $x = 0$ and $x=1$ through the points $A,B$ in the chart $\{ u_2 \neq 0\}$ satisfy:
    \begin{equation*}
        \begin{cases}
           u_4^2-4 u_2 u_5 =0,\\
           u_1^2 - 4 u_0 u_2 - 4 u_2 u_3 + 2 u_1 u_4 + u_4^2 - 4 u_2 u_5 = 0,\\
           u_5=0,\\
           u_0 + B_2 u_1 +B_2^2 u_2 + u_3 + B_2 u_4 + u_5 = 0,\\
           u_2 -1 = 0.
        \end{cases}
    \end{equation*}
    Let $V_{A,B,c,d}$ denote the zero set of this system and let  $(u_0,u_1,u_2,u_3,u_4,u_5)$ be an element of this family.    The line through $A$ and $B$ cannot be parallel to $x=0$ and hence there are no non-conic solutions.
    The third and fifth equations immediately imply that $u_5 = 0$ and $u_2 = 1$.
    After substituting these into the first equation gives $u_4^2 = 0$, that is, $u_4 = 0$.
    Substituting all of these in turn to the second and fourth equations gives the equations $u_1^2 -4u_0 -4u_3 = 0$ and $u_0 + B_2 u_1 + u_3 + B_2^2 = 0$. 
    Solving the first for $u_3$ and substituting to the second gives an equation only in the variable $u_1$,  which factors as $\frac{1}{4}(u_1 + 2 B_2)^2 = 0$.
    Hence, $u_1 = -2 B_2$.
    With these, the fourth equation becomes $-4u_0 -4u_3 + 4B_2^2 = 0$, that is, $-u_0 - u_3 + B_2^2 = 0$.
    Combining these observations, every element of $V_{A,B,c,d}$ is of the form $(u_0, -2 B_2, 1, B_2^2-u_0, 0, 0)$.
    Conversely, any point of this form satisfies all defining equations of the family $V_{A,B,c,d}$.
    Hence, the family $V_{A,B,c,d}$ is the line $\{(u_0, -2 B_2, 1, B_2^2-u_0, 0, 0) \mid u_0 \in \mathbb{R}\}$.
    In particular, $V_{A,B,c,d}$ is path-connected.

    Lastly, we show that the family $V_{A,B,c,d}$ contains a degenerate conic, namely the double line through the points $A$ and $B$.
    Let $l_{AB}$ be the line through $A$ and $B$ given by $B_2x-y = 0$.
    Then the corresponding double line is in $V_{A,B,c,d}$.
    This double line is a degenerate conic in $V_{A,B,c,d}$, and since $u_2 \neq 0$ for $E$ as well, $E$ is also contained in the same family, whose single path-connected component contains both the ellipse and this non-ellipse conic.     
\end{proof}

Influenced by the results of~\cite{kubjas2018positive}, one might ask whether it is enough to consider ellipses that have three incidences and two tangencies with the polygons $P$ and $Q$, respectively, or two incidences and three tangencies with the polygons $P$ and $Q$, respectively. We do not know the answer to this question.

\section{Semialgebraic description}\label{sec:semialgebraic_description}

In this section we give a semialgebraic description for the set $\mathcal{M}_{3,2}^{p \times q}$ of nonnegative matrices $p \times q$ matrices of rank at most $3$ and psd rank at most $2$.

\begin{theorem}\label{thm:semialgebraic-description}
    Let $M \in \R_{\geq 0}^{p \times q}$ be a nonnegative matrix with no zero rows or columns satisfying $M \mathbf{1} = \mathbf{1}$.
    Then $M \in \mathcal{M}_{3,2}^{p \times q}$ if and only if either $\rank(M)< 3$, or $\rank(M) = 3$ and, for a rank-size factorization $M = AB$ of $M$ where the rows of $A$ are of the form $(a_i,1)^T$, $a_i \in \R^2$, and columns of $B$ are $b_j \in \R^3$, there exist distinct $i_1, \dots, i_k \in [p]$, $j_{1},\dots,j_{5-k} \in [q]$ for $0 \leq k \leq 5$ and $u_0, \dots , u_5 \in \R$ defining a conic with the equation $U(x,y)=u_0x^2+u_1xy+u_2y^2+u_3x+u_4y+u_5$ and associated matrix $S$ such that 
    \begin{equation}\label{eqn:semialgebraic-description}
        \begin{cases}
             U(a_{(i,1)}, a_{(i,2)}) = 0 \text{ for } i\in \{i_1,\dots,i_k\},\\
             T(U,b_j)=0 \text{ for } j \in \{j_1,\dots, j_{5-k}\} ,\\
             U(a_{(i,1)}, a_{(i,2)})\leq 0 \text{ for } i\in [p] \backslash \{i_1,\dots,i_k\},\\
             T(U,b_j) \leq 0 \text{ for } j \in [q] \backslash \{j_1,\dots, j_{5-k}\},  \\
             u_0 > 0,\\
             4 u_0 u_2 -u_1^2 > 0, \text{ and} \\
             (u_0 +u_2) \det(S) < 0.
        \end{cases}
    \end{equation}
\end{theorem}
    
\begin{proof}
    By~\eqref{eqn:real-ellipse-conditions}, the last two conditions in \eqref{eqn:semialgebraic-description} assert that the conic defined by $U(x,y)$ is a real, nondegenerate ellipse.
    Since $U(x,y) = 0$ and $-U(x,y) = 0$ define the same conic, replacing the coefficients $u_0,\dots, u_5$ by $-u_0, \dots, -u_5$ does not change the described conic.
    Since $u_0 \neq 0$ by~\Cref{rmk:real_ellipse_nonzero_coefficients},
    we may assume without loss of generality that $u_0 >0$.
    The equation $U(a_{(i,1)}, a_{(i,2)})=0$ expresses that the vertex $a_i = (a_{(i,1)}, a_{(i,2)})$ of $P$ is on the conic $U$.
    By~\Cref{cor:point-inside-ellipse}, with the convention $u_0>0$, the inequality $U(a_{(i,1)}, a_{(i,2)}) < 0$ expresses that the vertex $a_i$ of $P$ lies on the interior of the conic.
    Similarly, the equation $T(U,b_j)=0$ expresses that conic $U$ is tangent to the line defined by the column $b_j$ of $B$, and by~\Cref{lem:inside-intersecting-ellipse}, a negative value $T(U,b_j) < 0$ expresses that the conic and the line do not intersect.

    First we prove the forward direction of the claim.
    Assume that $M \in \mathcal{M}_{3,2}^{p \times q}$.
    Then either $\rank(M) <3$ or $\rank(M) =3$. Assume that $\rank(M) = 3$ and let $M=AB$ be a size-$3$ factorization of $M$ where the rows of $A$ are of the form $(a_i,1)^T$ with $a_i \in \R^2$.
    Let $P$ and $Q$ be polygons as in \Cref{def:nested-polytopes-of-M}.
    As $\rank_{psd}(M) = 2$ it follows from~\Cref{prop:psd-rank-2-iff-ellipse} that there exists an ellipse $E$ such that $P \subset \operatorname{conv}(E) \subset Q$.
    By~\Cref{thm:main} there exists an ellipse $E'$ such that $P \subset \operatorname{conv}(E') \subset Q$ and $E'$ has at least five incidence and tangency relations with the polygons $P$ and $Q$.
    As $E'$ is a real, nondegenerate ellipse, it is defined by $U(x,y)$ such that the  coefficients $u_0, \dots, u_5 \in \R$ (multiplying with $-1$, if necessary) satisfy the last three conditions in \eqref{eqn:semialgebraic-description}.
    This means that there are vertices $a_{i_1}, \dots, a_{i_k}$, $0 \leq k \leq 5$, of $P$ and edges of $Q$ defined by $b_{j_1}, \dots, b_{j_{5-k}}$  such that the first and second equations in \eqref{eqn:semialgebraic-description} are satisfied.
    As $P \subset \operatorname{conv}(E')$, it follows that for the rest of the vertices $a_i$ of $P$ are contained inside or on the boundary of $E'$. 
    In particular, they satisfy the third inequality of \eqref{eqn:semialgebraic-description}.
    Similarly, since $\operatorname{conv}(E') \subset Q$, it follows that the rest of the edges $b_j$ of $Q$ intersect the ellipse in at most one point, that is, they satisfy the fourth inequality of \eqref{eqn:semialgebraic-description}.
    This proves the first direction of the claim.

    Next we prove the backward direction of the claim. Assume $M \in \R_{\geq 0}^{p \times q}$. 
    If $\rank(M) < 3$ then $\rank_{psd}(M) \leq 2$ by \Cref{prop:psd_rank_for_matrices_of_rank_one_and_two}, so $M \in \mathcal{M}_{3,2}^{p \times q}$.
    If $\rank(M) = 3$, consider a rank-size factorization $M=AB$ of $M$ where the rows of $A$ are of the form $(a_i,1)^T$, and let $u_0, \dots, u_5 \in \R$, $i_1, \dots, i_k \in [p]$, and $j_{1},\dots,j_{5-k} \in [q]$, for $0 \leq k \leq 5$, be such that the system \eqref{eqn:semialgebraic-description} is satisfied.
    The last two equations imply that $u_0, \dots, u_5$ are the coefficients of a real, nondegenerate ellipse $E$ defined by $U(x,y)$. 
    Since $u_0> 0$, the first and third conditions in \eqref{eqn:semialgebraic-description} imply that the vertices $a_{i_1}, \dots , a_{i_k}$ are on the boundary of the ellipse and the vertices $a_{i}$, $i\in [p] \backslash \{i_1,\dots,i_k\}$, are contained inside or on the boundary of the ellipse.
    In particular, the ellipse has  at least $k$ incidences with the polygon $P$, and, furthermore, $P \subset \operatorname{conv}(E)$.
    
    Similarly, the second and fourth conditions in \eqref{eqn:semialgebraic-description} express that the edges $b_{j_1}, \dots, b_{j_{5-k}}$ are tangent to the ellipse and the edges $b_j$ with $j \in [q]\backslash \{j_1, \dots, j_{5-k}\}$ intersect the ellipse in at most one point.
    In particular, the ellipse has at least $5-k$ tangencies with the polygon $Q$. By~\Cref{lemma:tact_invariant_polygon_containment}, $\operatorname{conv}(E) \subset Q$.
    Hence, $E$ is an ellipse such that $P \subset \operatorname{conv}(E) \subset Q$, so $M \in \mathcal{M}_{3,2}^{p \times q}$ by~\Cref{prop:psd-rank-2-iff-ellipse}.
    This proves the claim.
\end{proof}

\begin{remark}
This semialgebraic description can be easily written in the entries of the matrix $M$ and the auxiliary variables $u_0, \dots, u_5$ which correspond to coefficients of an ellipse even though the theorem is stated in terms of a rank-$3$ factorization $M=AB$ of $M$ when $\rank(M) = 3$.
    Suppose that $\rank(M) = 3$ and $M \mathbf{1}=\mathbf{1}$.
    The factorization can be constructed from $M$ by setting the columns of $A$ to be $M_k, M_l$ and $\mathbf{1}$ where $M_k$ and $M_l$ are columns of $M$ such that $A$ has rank three.
    This is possible because by assumption $\mathbf{1}$ is in the column space of $M$.
    Let $A^-$ be a left inverse of $A$.
    Then, we can set $B = A^-M \in \R^{3 \times q}$.
    The $j$th column of $B$ is $B_j = A^-M_j$ so that $M=AB$. 
\end{remark}

One can solve the first two conditions in \eqref{eqn:semialgebraic-description} for conics $U(x,y)$ that satisfy the given five incidences and tangencies, e.g., symbolically in~\texttt{Wolfram} or using numerical algebraic geometry software. For $M \in \R_{\geq 0}^{p \times q}$, there are at most $\sum_{k=0}^5 \binom{p}{k} \binom{q}{5-k}=\binom{p+q}{5}$ such systems. In practice, not all $a_i$ might give vertices of $P$ and not all $b_j$ might give edges of $Q$, so they can be removed from~\eqref{eqn:semialgebraic-description}, which reduces the number of systems considered.

The problem of how many reduced complex conics (i.e., all conics except double lines) pass through $k$ points and are tangent to $5-k$ lines in general position is investigated in \cite{bashelor2008enumerative}.
By~\cite[Table 2]{bashelor2008enumerative}, if $k \in \{0,5\}$, then there is one reduced conic, if $k \in \{1,4\}$, then there are two reduced conics, and if $k \in \{2,3\}$, then there are four reduced conics. The exact conditions for points and lines to be in general position are described in~\cite[Section 3]{bashelor2008enumerative}. In the case of three points and two lines,~\cite[Section 3.4]{bashelor2008enumerative} states that the points and lines are in general position unless the three points are collinear or a point lies on a line. In~\Cref{ex:inf-rigid-example}, we consider an example were these conditions are satisfied but the intersection point of the two lines is collinear with two of the points and we get that one of the reduced conics is a solution with multiplicity two. We think that this condition should be included among the configurations that are in general position. Nevertheless, one still gets four reduced conics when counted with multiplicities, so it does not change the counts above. In our setting, where points are vertices of $P$ and lines come from the edges of $Q$, unless a matrix $M$
 has a zero, we expect to get the number of solutions corresponding to reduced conics as stated above. 

Next, we consider examples of $3 \times 3$ matrices.
In these cases the polygons $P$ and $Q$ are triangles.

\begin{example}\label{ex:inf-rigid-example}
    The matrix in this example is a scalar multiple of the first example in \cite[Remark 6]{dawson2024uniquenesssize2positivesemidefinite}.
    Consider the following rank-$3$ matrix $M$ and its factorization $M = AB$  
    \begin{equation*}
        M = \begin{pmatrix}
            1/6 & 1/6 & 2/3 \\
            16/41 & 16/41 & 9/41\\
            4/41 & 36/41& 1/41
        \end{pmatrix} = \begin{pmatrix}
            2/3 & 0 & 1 \\
            16/41 & 16/41 & 1 \\
            0 & 16/41 & 1
        \end{pmatrix} \begin{pmatrix}
            3/4 & -5/4 & 1/2 \\
            53/48 & -5/16 & -19/24\\
            -1/3 & 1 & 1/3
        \end{pmatrix}.
    \end{equation*}
    The polygons associated to $M$ are $P = \operatorname{conv}\left\{(2/3 , 0) , (16/41 , 16/41), (0,16/41) \right\}$ and $Q = \{(x,y) \in \R^2 \mid -1/3 + 3 x/4 + 53 y/48 \geq 0, 1 - 5 x/4 - 5 y/16 \geq  0, 1/3 + x/2 - 19 y/24 \geq 0 \}$.
    We apply \Cref{thm:semialgebraic-description} to show that $\rank_{psd}(M) = 2$.
    By~\Cref{thm:semialgebraic-description} there are six possible configurations of five incidence and tangency relations between the conic and the polygons $P$ and $Q$: ellipses that pass through all three vertices of $P$ and are tangent to at least two edges of $Q$ (three configurations), or ellipses that pass through at least two vertices of $P$ and are tangent to all the edges of $Q$ (three configurations).

    We illustrate how to find an ellipse nested between $P$ and $Q$ through all vertices of $P$ and tangent to the edges of $Q$ contained in the lines $-1/3 + 3 x/4 + 53 y/48 = 0$ and $1 - 5 x/4 - 5 y/16 = 0$.
    This configuration corresponds to the following system: 
    \begin{equation}\label{eqn:rigid-example}
        \begin{cases}
            4 u_0/9 + 2 u_3/3 + u_5 = 0, \\
            256 u_0/1681 + 256 u_1/1681 + 256 u_2/1681 + 16 u_3/41 + 16 u_4/41 + u_5 = 0, \\
            256 u_2/1681 + 16 u_4/41 + u_5 = 0,\\
            
            256 u_1^2 + (53 u_3 - 36 u_4)^2 - 576 u_2 (4 u_3 + 9 u_5) + 16 u_1 (106 u_3 + 72 u_4 + 477 u_5) - 4 u_0 (256 u_2 + 848 u_4 + 2809 u_5) = 0, \\
            u_1^2 + 5/16 u_1 (2 u_3 + 8 u_4 + 5 u_5) + 1/256 (25 (u_3 - 4 u_4)^2 - 320 u_2 (4 u_3 + 5 u_5) - 4 u_0 (256 u_2 + 80 u_4 + 25 u_5))=0, \\
            64 u_1^2 + (19 u_3 + 12 u_4)^2 + 16 u_1 (19 u_3 - 12 u_4 - 57 u_5) + 192 u_2 (2 u_3 - 3 u_5) - 4 u_0 (64 u_2 + 152 u_4 + 361 u_5) \leq  0,\\
            u_0>0,\\
            4 u_0 u_2 -u_1^2 > 0, \\
            (u_0 +u_2) \det(S) < 0.
        \end{cases}
    \end{equation}
    The first three equations correspond to the equations $U(a_{(i,1)},a_{(i,2)})$ for the three vertices $p_i=(a_{(i,1)},a_{(i,2)}) \in \R^2$ of $P$.
    The fourth, fifth, and sixth equations are the tact invariants $T(U,b_j)$ of the three edges $b_j$ of $Q$, where the fourth and sixth equations are scaled by the positive scalars $2304$ and $576$, respectively, to clear the denominators.
    Solving the first five equations via~\texttt{Wolfram} gives three distinct solutions in $\mathbb{P}^5$. One solution has multiplicity $2$, so there are four solutions when counted with multiplicities matching the number in~\cite[Table 2]{bashelor2008enumerative}. 
    The corresponding conics are (we have chosen coordinates such that $u_0 > 0$):
        \begin{align*}
        & \frac{6 x^2}{7}+\frac{191 x y}{112} + \frac{697 y^2}{896} - x - \frac{29 y}{28} + \frac{2}{7}  = 0, \\
        & \frac{24 x^2}{35}  + \frac{1051 x y}{560} + \frac{17179 y^2}{13440} -x - \frac{599 y}{420} + \frac{38}{105} = 0, \text{ and } \\
        & \frac{3 x^2}{4} + \frac{29 x y}{16} + \frac{1681 y^2}{768} -x - \frac{41 y}{24} + \frac{1}{3}  = 0.
    \end{align*}
    The first solution is the degenerate conic consisting of a pair of intersecting lines whose intersection point lies at a vertex of $Q$.
    At this solution, the Jacobian of the system is not full-rank, and thus this is a solution with multiplicity $2$.
    The pair consists of the line through the vertices $(2/3 , 0)$ and $(16/41 , 16/41)$ of $P$ and the line through a vertex of $Q$ and $(0,16/41)$.
    The second solution is a hyperbola.
    The last solution satisfies the final two conditions of \eqref{eqn:rigid-example}, so it is a real ellipse: substituting the corresponding $U(x,y)$ to the sixth inequality in \eqref{eqn:rigid-example} yields $T(U,b_3) = 0$.
    Hence, by~\Cref{thm:semialgebraic-description} the matrix has psd rank $2$. 

    Repeating this analysis for the remaining five configurations, one finds in each case four solutions of conics, up to multiplicity and not counting double lines,  to the incidence and tangency equations of the system.
    Out of those exactly one is an ellipse nested between $P$ and $Q$.
    Moreover, the same ellipse occurs in each case.
    This is not surprising as the matrix $M$ has a unique, up to a certain group action of $\operatorname{GL}(2)$, psd factorization by \cite[Theorem 8]{dawson2024uniquenesssize2positivesemidefinite} which is equivalent to having a unique ellipse nested between, as each ellipse corresponds to a factorization.
    By~\cite[Proposition 3.8]{kubjas2018positive}, the matrix is on the boundary of $\mathcal{M}_{3,2}^{3 \times 3}$.
    The polygons $P$ and $Q$ and the ellipse $\frac{3 x^2}{4} + \frac{29 x y}{16} + \frac{1681 y^2}{768} -x - \frac{41 y}{24} + \frac{1}{3}  = 0$ are depicted in \Cref{fig:inf-rigid-example}.
         \begin{figure}[htbp!] 
         \centering
         \includegraphics[width=0.4\linewidth]{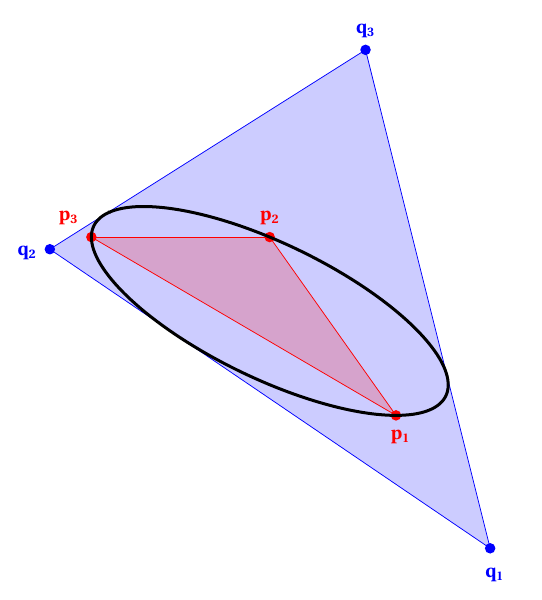}
         \caption{The polygons $P$ and $Q$ and the unique ellipse nested between them from~\Cref{ex:inf-rigid-example}.}
         \label{fig:inf-rigid-example}
     \end{figure}
\end{example}

\begin{example}
    In this example we present a matrix which does not have a unique ellipse nested between the associated polygons $P$ and $Q$.
    Let 
    \begin{equation*}
        M = \begin{pmatrix}
            1/2 & 3/8 & 1/8 \\
            3/10 & 1/40 & 27/40 \\
            1/8 & 3/8 & 1/2
        \end{pmatrix} = \begin{pmatrix}
            1/2 & 3/8 & 1\\
            3/10 & 1/40 & 1 \\
            1/8 & 3/8 & 1
        \end{pmatrix} \begin{pmatrix}
            1 & 0 & -1\\
            0 & 1 & -1\\
            0 & 0 & 1
        \end{pmatrix}.
    \end{equation*}
    The polygons associated to $M$ are $P = \operatorname{conv}\left\{(1/2, 3/8),(3/10, 1/40), (1/8, 3/8)\right\}$ and $Q = \{(x,y) \in \R^2 \mid (x,y,1)^T B \geq 0 \}$.
    By~\Cref{thm:semialgebraic-description} there are six possible configurations of five incidences and tangencies.
    Examining each configuration, we find four distinct solutions (not counting the double lines) out of which at least one is a nested ellipse in each case.
    In fact, we find several distinct ellipses that are nested between $P$ and $Q$:
    the ellipse given by $\frac{784}{25}x^2 -\frac{112}{5}xy + \frac{1024}{25}y^2 - \frac{56}{5}x-\frac{64}{5}y + 1 = 0$, which passes through all vertices of $P$ and is tangent to the edges of $Q$ contained in the lines $x=0$ and $y=0$, is nested between the polygons $P$ and $Q$.
    Similarly, the ellipse given by $\frac{4}{3}x^2  + \frac{4}{9}xy + \frac{4}{9}y^2  - x - \frac{4}{9}y + \frac{3}{16} = 0$, which passes through all vertices of $P$ and is tangent to the edges of $Q$ contained in the lines $y=0$ and $1-x-y=0$, is nested between the polygons $P$ and $Q$.
\end{example}

 \begin{example}
     In this example we present an example of a rank-$3$ matrix whose psd rank  is $3$.
     Let 
     \begin{equation*}
         M = \begin{pmatrix}
         13/15 & 1/15 & 1/15 \\
         1/15 & 13/15 & 1/15 \\
         1/15 & 1/15 & 13/15
     \end{pmatrix} = \begin{pmatrix}
         0 & 11/15 & 1 \\
         8/5 & -13/15 & 1 \\
         -8/5 & -13/15 & 1
     \end{pmatrix} \begin{pmatrix}
         0 & 1/4 & -1/4 \\
         1/2 & -1/4 & -1/4\\
         1/2 & 1/4 & 1/4
     \end{pmatrix}.
     \end{equation*}
     Then $P = \operatorname{conv}\left\{(0,11/15), (8/5, -13/15), (-8/5, -13/15)\right\}$ and $Q = \{(x,y) \in \R^2 \mid (x,y,1)B \geq 0 \}$.

     We use \Cref{thm:semialgebraic-description} to show that $\rank_{psd}(M) > 2$.
     There are again six configurations of five incidences and tangencies between the conic and the polygons. 
     In the first three configurations, where the conic passes through all the vertices of $P$ and is tangent to two edges of $Q$, all four solutions of the incidences and tangencies are hyperbolas.
     In the remaining three configurations, where the conic passes through two vertices of $P$ and is tangent to all edges of $Q$, all four solutions of the incidences and tangencies corresponding to reduced conics are real ellipses, but none of them satisfy the inequality $U(a_{(i,1)},a_{(i,2)}) \leq 0$ for the remaining vertex $p_i = (a_{(i,1)},a_{(i,2)})$ of $P$, that is, the vertex $p_i$ is not contained in the convex hull of the ellipse.
     Hence, $\rank_{psd}(M) > 2$ by \Cref{thm:semialgebraic-description}.
     Since $\rank(M) = 3$, it follows from \Cref{lem:relations-between-ranks} that $\rank_{psd}(M) = 3$.
 \end{example}

\noindent
\textbf{Acknowledgements.} We thank Jan Draisma for helpful discussions.
For 2026, the work of LM is supported by the Vilho, Yrjö and Kalle Väisälä Foundation of the Finnish Academy of Science and Letters.

\noindent \textbf{AI declaration.} The concept and all statements of the paper are by us. ChatGPT, Claude, and Proton LumoAI were used as tools in the process of writing this paper. We asked for ideas for proving~\Cref{thm:main}, which led us to the understanding that we have to study the geometry of conics with given incidence and tangent relations. The analysis was then done by us, with some help with details from AI models. We also used AI models for polishing text and converting mathematical objects from one input format to another.  Claude was used to generate two referee reports for the paper and we addressed many comments in these reports. We have checked everything and take full responsibility for the contents of the article.

\bibliographystyle{alpha}
\bibliography{references}

{\small{\noindent {\bf Authors' addresses}}

\medskip
\noindent{Kaie Kubjas, Aalto University}\hfill {\tt kaie.kubjas@aalto.fi}
\\ 
\noindent{Lilja Metsälampi, Aalto University}\hfill {\tt lilja.metsalampi@aalto.fi}
\\ 
\noindent{Milo Orlich, Aalto University}\hfill {\tt milo.orlich@aalto.fi}
\\ 
\end{document}